\documentclass[12pt]{article}
\usepackage{amsfonts}
\title{Lie Superalgebras, Clifford Algebras, Induced Modules and Nilpotent Orbits}
\author{Ian M. Musson\thanks{partially supported by NSF grant  DMS-0099923.}\\Department of Mathematical Sciences\\
University of Wisconsin-Milwaukee\\ email: {\tt
musson@csd.uwm.edu}}
\date{}
\newfont{\eufm}{eufm10 scaled\magstep1}
\newcommand{\FRAK}[1]{\mbox{\eufm#1}}
\font\twelveeufm=eufm10 scaled\magstep1 \font\teneufm=eufm10
\font\nineeufm=eufm9  \font\seveneufm=eufm7
\font\sixeufm=eufm6 \font\fiveeufm=eufm5
\newfam\eufmfam
\textfont\eufmfam=\twelveeufm \scriptfont\eufmfam=\nineeufm
\scriptscriptfont\eufmfam=\sixeufm \textfont\eufmfam=\teneufm
\scriptfont\eufmfam=\seveneufm
\scriptscriptfont\eufmfam\fiveeufm

\newtheorem{theorem}{Theorem}[section]
\newtheorem{lemma}[theorem]{Lemma}
\newtheorem{corollary}[theorem]{Corollary}

\renewcommand{\subsection}[1]{\refstepcounter{subsection}{~}\newline\textbf{\thesubsection.
#1}\quad}

\begin{document}
\date{ \today }
\maketitle
\begin{abstract}
Let $\FRAK{g}$ be a classical simple Lie superalgebra. To every
nilpotent orbit $\cal O$ in $\FRAK{g}_0$ we associate a Clifford
algebra over the field of rational functions on $\cal O$. We find
the rank, $k(\cal O)$ of the bilinear form defining this Clifford
algebra, and deduce a lower bound on the multiplicity of a
$U(\FRAK{g})$-module with $\cal O$ or an orbital subvariety of
$\cal O$ as associated variety.  In some cases we obtain modules
where the lower bound on multiplicity is attained using parabolic
induction.  The invariant $k(\cal O)$ is in many cases, equal to
the odd dimension of the orbit $G\cdot\cal O$ where $G$ is a Lie
supergroup with Lie superalgebra ${\mathfrak g.}$
\end{abstract}
\section  {\bf{\Large Introduction}}

Completely prime primitive ideals play a central role in  the
study of the enveloping algebra of a semisimple Lie algebra. For
example they are important in the determination of the scale
factor in Goldie rank polynomials, and they are related to unitary
representations, see [J3] for more details. On the other hand if
$\FRAK{g}$ is a classical simple Lie superalgebra, there are very
few completely
prime ideals in  $U(\FRAK{g}),$ see [M3, Lemma 1].\\

The results of this paper suggest that it may still be of interest
to study primitive ideals of low Goldie rank in  $U(\FRAK{g}),$
and their module theoretic analog, modules of low multiplicity.\\

To initiate this study we associate  to any prime ideal $q$ of
$S(\FRAK{g}_0),$ a Clifford algebra $C_q$ over the field of
fractions of $S(\FRAK{g}_0)/q.$ Let $k(q)$ be the rank of the
bilinear form defining this Clifford algebra.  Given a finitely
generated module $M,$ we use some filtered-graded machinery along
with an elementary result about Clifford algebras to obtain a
lower bound on the multiplicity of $M$ in terms of $k(q),$ see
Lemmas \ref{CA} and \ref{5.1}.

When $\FRAK{g}_0$ is reductive and $P$ is a primitive ideal in
$U(\FRAK{g}_0)$ the subvariety of $\FRAK{g}_0$ defined by $gr P$
is the closure of a nilpotent orbit, [BB], [J1].  For this reason
the most interesting primes in $S(\FRAK{g}_0)$ are those defining
nilpotent orbits or their orbital subvarieties.
 If $\FRAK{g}$ is classical simple and $q$ is a prime
 ideal of  $S(\FRAK{g}_0)$ defining a nilpotent orbit we give
 a formula for $k(q)$ in terms of a partition (or
 partitions) associated to the nilpotent orbit.\\
\indent This work motivates the search for highest weight modules
with given associated variety and low multiplicity.  For $\FRAK{g}
= g\ell(m,n), \; s \ell(m,n)$ or $Q(n)$ we explain how to find
examples of such modules using induction from parabolic
subalgebras.  For a precise statement, see Lemmas \ref{5.4},
\ref{5.5} and Theorem \ref{5.6}.  We also investigate the
primitive ideals that arise as annihilators of these modules and
the structure of the corresponding primitive factor algebras. We
remark that the orbital varieties which occur in our examples have
the simplest possible type, namely they are all linear
subvarieties of the nilpotent orbit.  One difficulty is that the
closest analog for semisimple Lie algebras of the problem
considered here is the quantization problem for orbital varieties,
which is unsolved, see [Be],[J3]. It is worth noting also that the
associated variety of a simple highest weight module is
irreducible for $sl(n)$, [Me]. This is not true in general
[J2],[T]. We plan to return to the issues raised here in a
subsequent paper.  In particular we shall show that the modules we
construct in this paper are quantizations of superorbital
varieties.

Additional motivation for the study of the invariants $k(q)$ comes
from supergeometry. Suppose that $\FRAK{g}$ is classical simple,
and that there is a nondegenerate even bilinear form on
$\FRAK{g}.$ If $x \in \FRAK{g}_0$, and $m_x$  is the corresponding
ideal of $S(\FRAK{g}_0)$ then $k(m_x)$ is equal to the dimension
of the centralizer of $x$ in $\FRAK{g}_1.$ If $G$ is a Lie
supergroup with Lie superalgebra $\FRAK{g}$, this allows us to
find the superdimension of the orbit $G\cdot x$, when $x$ is
nilpotent.

 \indent This paper is organized as follows. After some
preliminaries in Section 2, we obtain our formulas for $k(q)$ in
Section 3.  Although this is done on a case-by-case basis, the
formulas in most cases depend on the same basic result (Lemma
\ref{Part}).  Furthermore the exceptional algebras $G(3)$ and
$F(4)$ can be treated using essentially the same method as the
orthosymplectic algebras.  In Section 5 we prove our main results
about parabolically induced modules. We prove a result (Theorem
\ref{5.2}) describing the structure of such modules as
$U(\FRAK{g}_0)$-modules.  This is used to derive analogs of
several results on induced modules and their annihilators from
[Ja2, Kapitel 15].  Several  of the results in this section (for
example Theorem 5.7 and Corollary 5.10) apply to the modules
$F(\mu)$ constructed by Serganova in section $3$ of [S2] for the
Lie superalgebras $g\ell(m,n).$ In Section 4 we give some
background on parabolic subalgebras needed in Section 5. Our
results on nilpotent orbits may be found in Section 6.
 Nilpotent orbits do not seem to have been
widely studied in the superalgebra case, see however [S1], so we
spend some time developing the background.\\

\indent I would like to thank Olivier Mathieu,  Vera Serganova and
Jeb Willenbring
for some useful discussions.\\
\\
\\
\section {\Large{\bf Preliminaries}} \label{test}
\subsection {\bf  Clifford Algebras.} {\label{CA}}
Let $\FRAK{g} = \FRAK{g}_0 \oplus \FRAK{g}_1$ be a finite
dimensional Lie superalgebra over $\mathbb{C}$.  The tensor
algebra $T(\FRAK{g})$ has a unique structure
 $T(\FRAK{g}) = \oplus_{n \geq 0}T^n(\FRAK{g})$ as a graded algebra such
 that $T^0(\FRAK{g}) = \mathbb{C}, T^1(\FRAK{g}) = \FRAK{g}_1$ and
 $T^2(\FRAK{g}) = \FRAK{g}_0 + \FRAK{g}_1 \otimes \FRAK{g}_1$.
 Set $T_n = \oplus_{m \leq
 n}T^m(\FRAK{g})$ and let $U_n$ be the image of $T_n$ in $U(\FRAK{g})$.
Then $\{ U_n \}$ is a filtration on $U(\FRAK{g})$ and we describe
the associated graded ring $S = gr U(\FRAK{g})$. Observe that $R =
S(\FRAK{g}_0)$ is a central subalgebra of $S$ and that the bracket
$[\;,\;]$ on $\FRAK{g}_1$ extends to an $R$-bilinear form on
$\FRAK{g}_1 \otimes R$.  The algebra $S$ is isomorphic to the
Clifford algebra of this bilinear form. If $v_1, \ldots, v_n$ is a
basis of $\FRAK{g}_1$ over $\mathbb{C}$ then the matrix of the
bilinear form with respect to this basis is $M(\FRAK{g}) = ([v_i,
v_j])$. We do not refer to the basis in the notation for this
matrix since we study only properties of the matrix which are
independent of the basis.

 \indent
We showed in [M2] that there is a homeomorphism
 \[ \pi : Spec R
\longrightarrow Gr Spec S.\] where $Gr Spec (\;\;)$ refers to the
space of $\mathbb{Z}_2$-graded prime ideals.  Let us recall the
details.  Fix $q \in Spec(R)$ and let
 $\overline{S} = S/Sq$ and $\mathcal{C} = \mathcal{C}(q)$, the set of regular
elements of $R/q$.  Then $F_q = Fract(R/q)$ is a central subfield
of the localization $T = \overline{S}_{\cal C}$.  Moreover the Lie
bracket on $\FRAK{g}_1$ extends to a symmetric $F_q$-bilinear form
on $\FRAK{g}_1 \otimes F_q$. It is easy to see that $T$ is the
Clifford algebra of this form over $F_q$.  The nilradical $N$ of
$T$ is generated by the radical of the bilinear form on
$\FRAK{g}_1 \otimes F_q$, and $T/N$ is the Clifford algebra of a
nonsingular bilinear form. Then $\pi(q)$ is the kernel of the
combined map
 \[ S = gr U(\FRAK{g}) \longrightarrow \overline{S}
 \longrightarrow T/N.  \]
 It follows that $\pi(q) =  \sqrt{Sq}$ where  $\sqrt{\;\;\;}$ denotes the radical of an
 ideal.
 For $p \in Gr Spec S, \pi^{-1}(p) = p
\cap R$. Note that if $\pi(q) = p$ we have inclusions of rings
\[ R/q \subseteq S/p \subseteq T/N.\]
Moreover $T/N$ is obtained from $S/p$ by inverting the nonzero
elements of $R/q$. Hence $S/p$ is an order in the Clifford algebra
$C_q = T/N$. Let $B_q$ be the bilinear form defining this Clifford
algebra, $\delta_q$ the determinant of $B_q$ and $k(q)$ the rank
of $B_q.$ Thus
\[ k(q) = \{ max \; m|\mbox{some} \; m \times m \;\mbox{minor of}\; M(\FRAK{g})\;
  \mbox{is nonzero mod}\, q\}. \]

 A prime ideal $q$ of $S(\FRAK{g}_0)$ is {\it homogeneous} if
$q = \oplus_{n \geq 0}(q \cap S^n(\FRAK{g}_0))$ where
$S(\FRAK{g}_0) = \oplus_{n \geq 0}S^n(\FRAK{g}_0)$ is the usual
grading.  All prime ideals $q$ of $S(\FRAK{g}_0)$ considered in
this paper will be homogeneous. If $q$ is homogeneous and $k(q)$
is odd then $\delta_q$ is a rational function of odd degree and
hence not a square in $F_q$. Therefore by [L, Theorems V.2.4 and
V.2.5] $C_q$ is a central simple algebra.  Hence $C_q \cong
M_{2^a}(D)$ for a division algebra $D$. Using the fact that
$\dim_{F_q} C_q = 2^{k(q)}$ it is easy to prove the following
result.
\begin{lemma} {\label{sz}}
Let $L$ be a simple $C_q$-module where $q$ is a homogeneous
prime ideal of $S(\FRAK{g}_0)$.  \\
 (a) If $k(q)$ is even then $C_q$ is a central simple algebra
 over $F_q$ and $\dim_{F_q} L \geq 2^{k(q)/2}$.  Equality holds if and only if
 $D = F_q$\\
 \\
 (b) If $k(q)$ is odd then $C_q$ is a central simple algebra
 over $F_q(\sqrt{\delta_q})$ and $\dim_{F_q}L \geq
 2^{(k(q)+1)/2}$.  Equality holds if and only if $D = F_q(\sqrt{\delta_q})$.
 \end{lemma}

 We denote the greatest
 integer less than or equal to $s$ by [s].  If $\dim_{F_q} L = 2^{[(k(q)+1)/2]}$,
 we say that $C_q$ is {\it split}.  \\
 \\
\subsection {Evaluation of  $M(\FRAK{g}).$}  \label{EM}
 Let $\FRAK{g}$ be
classical simple.  Since $\FRAK{g}_0$ is reductive there is a
nondegenerate invariant bilinear form on $\FRAK{g}_0 $.  This
allows us to identify $\FRAK{g}_0$ with $\FRAK{g}^*_0$ and
 thus to view elements of $S(\FRAK{g}_0)$ as functions on
 $\FRAK{g}_0$.
  If ${\cal O}  \subseteq \FRAK{g}_0,$ and the ideal $q$ of
  functions in $S(\FRAK{g}_0)$ which vanish on  ${\cal O}$ is prime,
  we often write $k(\cal O)$ in place of $k(q)$.
It is convenient to set $\ell(q) = [(k(q) + 1)/2]$ and $\ell({\cal
O}) = [(k({\cal O}) + 1)/2]$.  We say that a closed subset $X$ of
$\FRAK{g}_0$ is {\it conical} if $x \in X$ implies that
$\mathbb{C}x \subseteq X.$  For example closures of nilpotent
orbits and their orbital subvarieties are conical.  If $X$ is a
product of conical subvarieties of the simple summands of
${\mathfrak g}_0,$ then the defining ideal of $X$ in $S({\mathfrak
g}_0)$ is independent of the choice of bilinear form, since any
two nondegenerate invariant forms on a simple Lie algebra are
proportional.
 Fix a nilpotent orbit $\cal O$, and suppose $q \in Spec
S(\FRAK{g}_0)$ is such that $V(q) = {\overline{\cal O}}$.  We want
to compute $k(q)$.  For $x \in {\cal O}$, let $M(x)$ be the
evaluation of $M(\FRAK{g})$ at $x$ and let $m_x$ be the maximal
ideal of $S(\FRAK{g}_0)$ corresponding to  $x$.  Since ${\cal O}$
is dense in $V(q)$ and the rank of $M(\FRAK{g})$ is constant on
${\cal O}$ we have
  \begin{equation}\label{new}
 k(q) = rank(M(x)) = k(m_x)  \quad \mbox{for all} \quad x \in {\cal   O}.
  \end{equation}
Hence if $X$ is an irreducible subvariety of $\mathcal{O}$ we have
$k(X) = k(\mathcal{O})$.
 \\
 \subsection {\bf Matrix notation.}  \label{MN}
We denote the $n \times n$ identity matrix by $I_n,$ and the
matrix with a $1$ in row $i$, column $j$ and zeroes elsewhere by
$e_{ij}.$  Let $\Upsilon_r$ be the $r \times r $ matrix with ones
on the antidiagonal and zeros elsewhere. We write  $M_{m,n}$ for
the vector space of $m \times n$ complex matrices. The transpose
of a matrix $A$ is denoted by $A^t.$ Since $M(\FRAK{g})$ is a
matrix over $S(\FRAK{g}_0),$ and $\FRAK{g}_0$ is often an algebra
of matrices, we need an "external" version of the matrices
$e_{ij}.$  For clarity, a matrix $A$ with entries in
$S(\FRAK{g}_0)$ will often be written in the form
\[ A = \sum_{i,j} a_{i,j} {\bf e_{i,j}} \]
meaning that $a_{i,j} \in S(\FRAK{g}_0)$ is the entry in row $i$
and column $j$ of $A.$

 Recall that if $A$ and $B$ are square matrices with rows and
columns indexed by $I,J$ respectively, then the Kronecker product
$A \otimes B$ has rows and columns indexed by $I \times J$, and
has entry in row $(i,k)$, column $(j, \ell)$ equal to $a_{ij}b_{k
\ell}$. To be more precise, we should also specify an ordering on
the rows and columns of  $A \otimes B.$  If $I \subseteq
\mathbb{Z}$ we give $I$ the ordering inherited from $\mathbb{Z}$.
If $I,J$ are ordered sets then unless otherwise stated we give $I
\times J$ the lexicographic order $<_{\ell ex}$ defined by
 \[ (i,j) <_{\ell ex} (k,\ell) \quad \mbox{if and only if} \quad i < k
 \quad\quad \mbox{or}\quad i = k \quad \mbox{and} \quad j < \ell . \]

We need a twisted version of the Kronecker product.  If $A$ and
$B$ are as above, we define  $A \widehat{\otimes} B$ to be the
matrix with rows indexed by $I \times J$ and columns indexed by $J
\times I$ such that the entry in row $(i,k)$, column $(\ell,j)$ is
equal to $a_{ij} b_{\ell k}$. Here we order $I \times J$ the
lexicographically  and order $J \times I$ so that $(j,i)$ precedes
$(\ell,k)$ if and only if $(i,j) <_{\ell ex} (k,\ell)$ .

The definition of $A \widehat{\otimes} B$ might seem unnatural at
first, but it is very convenient for the computation of
$M(\FRAK{g})$ when $\FRAK{g} = g\ell(m,n)$.  Note that if we
relabel column $(\ell, j)$ of $B$ as column $(j, \ell)$, the rows
and columns of $A \widehat{\otimes} B$ are then both indexed by $I
\times J$ ordered lexicographically.  It follows that $A
\widehat{\otimes} B = A \otimes B^t$.
\\
 \\
\subsection {\bf  Partitions.} {\label{Part}}
  If $\mu = (\mu_1 \geq \mu_2
\geq \ldots)$ is a partition of $m$ we denote the nilpotent matrix
with Jordan blocks of size $\mu_1, \mu_2, \ldots$  by  $J_\mu$.
The dual partition $\mu'$ of $m$ is defined by
\[ \mu'_i = |\{ j| \mu_j \geq i \}|  \]
for all $i$. We set $\mu_i = 0$ for all $i  > \mu'_1.$
The set of all partitions of $m$ is denoted ${\bf P}(m)$.\\
 \begin{lemma} For $\mu \in {\bf P}(m)$ and $\nu \in  {\bf  P}(n)$ we have
\[rank (J_\mu \otimes I_n + I_m \otimes J_\nu) =  mn - \sum_{i \geq 1} \mu'_i \nu'_i  . \]
 \end{lemma}
\noindent
 {\bf Proof.}
For $a \geq 1$, let $L(a)$ be the simple $s\ell(2)$-module of
dimension $a$.  If $e = \left[ \begin{array}{cc} 0 & 1\\0 & 0
\end{array}\right]$ we can choose bases for the modules
$L(\mu_i)$ and $L(\nu_i) $ such that $E= J_\mu \otimes I_n +
I_m\otimes J_\nu$ is the matrix representing the action of $e$ on
\[ \oplus_{i \geq 1} L(\mu_i) \otimes \oplus_{i \geq 1} L(\nu_i) . \]
To compute {\it rank} $E$ note that $L(a) \otimes L(b)$ is the
direct sum of $\min (a,b)$ simple modules, and the rank of $e$
acting on $L(a)$ is $a-1$.  This implies\\
 \begin{eqnarray*} \label{e2}
\quad\quad\quad\quad\quad\quad {\it rank} E = mn - \sum_{j,k}
\min(\mu_j,\nu_k).
  \end{eqnarray*}
  Now set
\[ A_i = \{ (j,k)| \min(\mu_j, \nu_k) = i \} \]
\[ B_i = \{ (j,k)| \min(\mu_j, \nu_k) \geq i \} .\]
Note that $|B_i| = \mu'_i \nu'_i$ and $|A_i| = |B_i| - |B_{i+1}|$.
Thus\\
 \begin{eqnarray*} \label{e3}
  \quad\quad\quad\quad\quad  \sum_{j,k}\min(\mu_j,\nu_k) = \sum_i
i|A_i| = \sum_i |B_i| = \sum_i \mu'_i \nu'_i .
  \end{eqnarray*}
 \\
{\bf Remark.} Since  $J_\mu \widehat{\otimes} I_n = J_\mu \otimes
I_n$ and  $I_m \widehat{\otimes} J_\nu = I_m \otimes J_\nu^t$ we
also have a formula for  $rank (J_\mu
\widehat{\otimes} I_n + I_m \widehat{\otimes} J_\nu).$\\

 \subsection {Dimension and Multiplicity.}  {\label{DM}}
  Let $N = \oplus_{m \geq 0} N(m)$ be a finitely generated graded $S(\FRAK{g}_0)$-module
 and set $N_n = \oplus^n_{m = 0} N(m)$.  For $n > > 0$ we
  have
  \[ \dim N_n = a_d \left( \begin{array}{cc} n\\d
  \end{array}\right) + a_{d-1}\left( \begin{array}{cc} n \\ d-1
  \end{array}\right) + \ldots + a_0 \]
  for suitable constants $a_0, \ldots, a_d$ with $a_d \neq 0$.
We set $d(N) = d$ and $e(N) = a_d.$ We filter $U(\FRAK{g})$ as in
Section \ref{CA} and denote associated graded ring by $gr
U(\FRAK{g})$. Let $M$ be a finitely generated $U(\FRAK{g})$-module
and
 equip $M$ with a good filtration $\{M_n\}_{n \geq 0}$.
 Since $N = gr M$ is finitely generated over $gr U(\FRAK{g})$
  and hence over $S(\FRAK{g}_0)$, the above remarks apply and we
  set $d(M) = d(N)$ and $e(M) = e(N).$
 It is not hard to show that $d(M)$ and $e(M)$ are independent of the
  good filtration and that $d(M)$ is the Gelfand-Kirillov
  dimension of $M$ calculated either as a $U(\FRAK{g})$-module or as
  a $U(\FRAK{g}_0)$-module. For details see [KL, Chapter 7].
If $M$ is finite dimensional, we have $d(M) = 0$ and $e(M) =
\dim_{\mathbb{C}}M$.

 In Section 5 we use the following  fact.  Suppose $q$ is a
homogeneous prime ideal of $S(\FRAK{g}_0)$ and $N$ a finitely
generated torsionfree graded module over $Z = S(\FRAK{g}_0)/q$. If
$F = Fract(Z)$ then $d(N) = d(Z)$ and $e(N) = e(Z) \cdot \dim_F
Z^{-1} N$.  This follows easily from [GW, Exercise 4L, Corollary
4.17 and Lemma 6.17].  If ${\cal {V}}$ is the closed subset of
$\FRAK{g}^*_0$ defined by $q$ we set $e({\cal{V}}) = e(Z)$.

A module $M$ is {\it homogeneous} (resp. {\it critical}) if for
any nonzero submodule $M'$ we have $d(M) = d(M')$ (resp. $d(M) =
d(M')$ and $e(M) = e(M'))$.
   \\
   \\
   \subsection {\bf  Induced Modules.} {\label{IM}}
    Let ${\FRAK{p}}$ be a subalgebra of the Lie superalgebra
    $\FRAK{g}$ and $N$ a finitely generated $U(\FRAK{p})$-module.
    We write $Ind^{\FRAK{g}}_{\FRAK{p}}\; N$ for the induced module
     $U(\FRAK{g}) \otimes_{U(\FRAK{p})}N$.
     \begin{lemma}
     Suppose $M = Ind^{\FRAK{g}}_{\FRAK{p}}\:N$ and set $c_i = \dim
     \FRAK{g}_i - \dim \FRAK{p}_i$ for $i = 0,1$. Then
     \[ d(M) = d(N) + c_0 \]
     and
     \[ e(M) = 2^{c_1} e(N) . \]
 \end{lemma}
{\bf Proof.} This is easily adapted from the proof of [Ja2, Lemma
8.9].
 \\
 \\
 \subsection {\bf  Affiliated Series of a Module.}  \label{AS}
 Let $N$ be a nonzero finitely generated module over a Noetherian
 ring $S.$
 An {\it affiliated submodule} of $N$ is a submodule of the
 form $ann_N(P)$ where $P$ is an ideal of $S$ maximal among the
 annihilators of nonzero submodules of $N$, see [GW] for background.  An {\it affiliated
 series} for $N$ is a series of submodules
  \[ 0 = N_0 \subset N_1 \subset \ldots \subset N_k = N \]
  such that each $N_i/N_{i-1}$ is an affiliated submodule of
  $N/N_{i-1}$.  The prime ideals $P_i = ann_S(N_i/N_{i-1})$ are
  called the {\it affiliated primes} of the series.

\subsection {\bf  Reductive Lie algebras.}{\label{rla}} For the remainder
of section 2,  $\FRAK{g}_0$ will be a reductive Lie algebra. Later
we use the notation established here when  $\FRAK{g}_0$ is the
even part of a classical simple Lie superalgebra. Let $\FRAK{n}_0
\oplus \FRAK{h}_0 \oplus \FRAK{n}^+_0$ be a triangular
decomposition of $\FRAK{g}_0$. So $\FRAK{h}_0$ is a Cartan
subalgebra and ${\bf b} = \FRAK{h}_0 \oplus \FRAK{n}^+_0$ a Borel
subalgebra of $\FRAK{g}_0$.  Let $G$ be the adjoint algebraic
group of $\FRAK{g}_0$.  If $\alpha$ is a root of $\FRAK{g}_0$ we
denote the corresponding root space by $\FRAK{g}^\alpha.$  There
is a unique element $h_{\alpha} \in
[\FRAK{g}^{\alpha},\FRAK{g}^{-\alpha}] $ such that $\alpha
(h_{\alpha}) = 2.$ For $\lambda \in \FRAK{h}_0^*$ we denote the
Verma module with highest weight $\lambda$ induced from ${\bf b}$
and its unique simple quotient by $M(\lambda)$ and $L(\lambda)$
respectively. We write $(\lambda,\alpha^{\vee})$ in place of
$\lambda(h_{\alpha}).$

\subsection {\bf Richardson Orbits.}{\label{rich}}
Let $\FRAK{p}_0$ be a parabolic subalgebra of $\FRAK{g}_0$ and
suppose that $\FRAK{p}_0 = \FRAK{l}_0 \oplus \FRAK{m}_0$ where
$\FRAK{m}_0$ is the nilradical of $\FRAK{p}_0$ and $\FRAK{l}_0$ is
a Levi factor. Then $G \FRAK{m}_0$ contains a unique dense orbit
called the {\it Richardson orbit induced from
$\FRAK{l}_0$ }.\\

\indent If $L$ is a finite dimensional $\FRAK{l}_0$-module and $M
= Ind^{\FRAK{g}_0}_{\FRAK{p}_0}\; L$ there are two prime ideals of
$S(\FRAK{g}_0)$ that we can associate to $M$.  The first of these
is $q' = \sqrt{ gr ann_{U(\FRAK{g}_0)}M}$ which is the defining
ideal of the Richardson orbit $\cal{O}$ induced from $\FRAK{l}_0$,
[Ja2, 17.15] . On the other hand we can equip $M$ with a good
filtration and consider $q = \sqrt{ ann_{S(\FRAK{g}_0)} gr M}.$
Then $q = S(\FRAK{g}_0)\FRAK{p}_0$ is the defining ideal of
$\FRAK{m}_0 \subset \overline{\cal{O}}$, [Ja2, 17.12 (4)]. We have
$2 dim ( \FRAK{m}_0) = dim(\cal{O})$. However $k(q) = k(q')$ since
$\cal{O} \; \cap$ $ \FRAK{m}_{0}$ is nonempty and by (\ref{new})
in section \ref{EM} $k(q)$ can be calculated by evaluating at any
point
of $\cal{O}.$\\
\\
\subsection{\bf Orbital Varieties.}{\label{OV}}
Let $\cal O$ be  a nilpotent orbit in $\FRAK{g}_0$.  The
irreducible components of ${\cal O} \cap \FRAK{n}_0^+$ are called
{\it orbital varieties attached to} ${\cal O}$. If ${\cal V}$ is
such an orbital variety we have $k({\cal O}) = k({\cal V})$ as
above.  For example if $\cal O$ is the Richardson orbit induced
from $\FRAK{l}_0$ and $\FRAK{m}_0$ is as in section \ref{rich}
then $\FRAK{m}_0$ is an orbital variety in $\cal O$.  In general
however Richardson orbits contain many other orbital varieties,
see [J3] for a recent survey.

\subsection {\bf The category O.} {\label{Cat}}
We denote by {\bf O} the category of  $U(\FRAK{g}_0)$-modules
defined in [Ja1, section 1.9].  For $M \in Ob \; {\bf O} $ we
write $[M]$ for the class of $M$ in the Grothendieck group $G({\bf
O})$ of ${\bf O}$.  The group  $G({\bf O})$ is free abelian on the
classes $[ L(\lambda) ]$ with $\lambda \in \FRAK{h}_0^*$ . For $M,
M' \in Ob \; {\bf O}$ we have $[M] = [M']$ if and only if $M$ and
$M'$ have the same character. We define a partial order $\leq$ on
$G({\bf O})$ by the rule $ \sum_{\lambda} a_{\lambda}[L(\lambda)]
\leq \sum_{\lambda} b_{\lambda}[L(\lambda)]$ if and only if
$a_{\lambda} \leq b_{\lambda}$ for all $\lambda \in \FRAK{h}_0^*.$

\section {\Large{\bf Dimension Formulas}}
\subsection {\label{3.1}}
 We describe the matrix $M(\FRAK{g})$ explicitly when $\FRAK{g} = g \ell(m,n)$. Let ${\bf I_1} = \{
1, \ldots, m\}, {\bf I_2} = \{ m+1, \ldots, m+n\}$, ${\bf I} =
{\bf I_1}  \cup {\bf I_2}$ and consider the following matrices
\[ N_1 = \sum_{i,j \in {\bf I_1} } e_{ij} {\bf e_{i,j}}, \quad N_2 = \sum_{k,\ell \in
{\bf I_2}} (e_{k,\ell}){\bf e_{k,\ell}} ,\]
 with entries in $\FRAK{g}_0$.\\
\\
\begin{lemma}
With a suitable choice of ordered basis for $\FRAK{g}_1,
M(\FRAK{g})$ has block matrix form
\[  \left[ \begin{array}{cc}
    0     & N \\
   N^{t} & 0
  \end{array} \right] \]
 where $N = N_1 \widehat{\otimes} I_n + I_m \widehat{\otimes} N_2.$
\end{lemma}
{\bf Proof.} Write $\FRAK{g}^+_1 = span\{e_{ik}|(i,k) \in {\bf
I_1} \times {\bf I_2} \}$, $\FRAK{g}^-_1 = span\{e_{\ell
j}|(\ell,j) \in {\bf I_2} \times {\bf I_1} \}$, so that
$\FRAK{g}_1 = \FRAK{g}^+_1 \oplus \FRAK{g}^-_1$.  The rows and
columns of $M(\FRAK{g})$ are indexed by ${\bf I_1} \times {\bf
I_2}$ ordered lexicographically followed by ${\bf I_2} \times {\bf
I_1}$ ordered so that $(j,i)$ precedes $(\ell,k)$ if and only if
$(i,j) <_{\ell ex} (k,\ell)$ .

The block matrix decomposition follows since $[\FRAK{g}^\pm_1,
\FRAK{g}^\pm_1] = 0$ and $M(\FRAK{g})$ is symmetric.  To compute
$N$ suppose $(i,k) \in {\bf I_1} \times {\bf I_2}$ and $(\ell, j)
\in {\bf I_2} \times {\bf I_1}$, then $[e_{ik}, e_{\ell j}] =
\delta_{k \ell} e_{ij} + \delta_{ij} e_{\ell k}$ and the result
follows.
 \\
 \\
\subsection {\label{3.2}}
Let $N_1 = \sum_{i,j \in {\bf I_1} } e_{ij} {\bf e_{i,j}}$
as above and $y = \sum_{k,\ell}y_{k\ell} e_{k\ell} \in g \ell
(m)$.  Using the bilinear form $(A,B) = trace(AB)$ to evaluate
$N_1$ at $y$ we have that
\[ N_1 (y) = (y_{ji}) \in g \ell (m) \]
is the $m \times m$ matrix with $i,j$ entry equal to $y_{ji}$.
Thus $N_1(y)$ has the same Jordan form as $y$.  Of course similar
remarks apply to the evaluation of $N_2.$

We denote the orbit of $(J_\mu, J_\nu)$ in $\FRAK{g}_0 = g\ell(m)
\times g\ell(n)$ by ${\cal O}_{\mu,\nu}$.\\

 \begin{theorem}
For $\mu \in {\bf  P}(m)$ and $\nu \in {\bf  P}(n)$ we have
\[ k({\cal O}_{\mu,\nu}) = 2(mn - \sum_i \mu'_i \nu'_i ) . \]
 \end{theorem}
  \noindent{\bf Proof.} This is immediate by Lemmas \ref{Part} and
 \ref{3.1}.
  \\
  \\
{\bf Remark.} If $\FRAK{g} = s\ell(m,n)$ then $\FRAK{g}$ has the
same odd part as $g\ell(m,n)$ and the matrix $M(\FRAK{g})$ can be
calculated using Lemma \ref{3.1}.  We can identify the nilpotent
orbits in  $\FRAK{g}_0$ with those in the even part of
$s\ell(m,n)$ and then Theorem  \ref{3.2} applies to $\FRAK{g}$  .
Similar remarks apply to the Lie superalgebra $ps\ell(n,n)$.\\
 \\

\subsection {\label{3new}}
If $V$ is a vector space we write  $\wedge^k V$ and $S^k V$ for
the $kth$ exterior and symmetric power of $V$ respectively. For
$v, w \in V$ we set $v \wedge w = 1/2 (v \otimes w - w \otimes v)
\in \wedge^2 V,\; v \circ w = 1/2(v \otimes w + w \otimes v) \in
S^2 V$. The following description of the orthosymplectic Lie
superalgebra algebra $osp(m,n)$ can be found in [K, 2.1.2]. Let
$V_1$ be an $m$-dimensional vector space with a nondegenerate
symmetric bilinear form $\psi_1$ and $V_2$ an $n$-dimensional
vector space with a nondegenerate skew-symmetric bilinear form
$\psi_2$.

 Then we can realize $ \FRAK{g} = osp(m,n)$ by setting
\[  \FRAK{g}_0 = \wedge^2 V_1 \oplus S^2 V_2, \;  \FRAK{g}_1 = V_1 \otimes V_2. \]
The action of $\wedge^2 V_1$ on $V_1$ is given by
\[ [a \wedge b, c] = \psi_1(a,c)b - \psi_1(b,c)a  . \]
Similarly $S^2 V_2$ acts on $V_2$ via
\[ [a\circ b,c] = \psi_2(a,c) b + \psi_2(b,c)a .\]
The bilinear forms $\psi_1$ and $\psi_2$ are invariant under these
actions, so $\wedge^2 V_1$ and $S^2 V_2$ identify with $so(m)$ and
$sp(n)$ respectively. The product  $ \FRAK{g}_1 \times  \FRAK{g}_1
\longrightarrow \FRAK{g}_0$ is given by
\[ [a \otimes c, b \otimes d] = \psi_1(a,b) (c \circ d) + \psi_2(c,d)(a \wedge b) .\]

\subsection {\label{3.3}}
The following lemma applies to the computation of the matrix
  $M(\FRAK{g})$ when $\FRAK{g}_0$ is not simple, $\FRAK{g}_1$
  is an irreducible $\FRAK{g}_0$-module and $\FRAK{g} $ is not
  isomorphic to  $\Gamma(\sigma_1, \sigma_2, \sigma_3).$
The discussion leading up to  [Sch, page 143, equation (5.9)]
allows us to make the following assumptions about the structure of
$\FRAK{g}.$ Firstly $\FRAK{g}_0 = \FRAK{g}^1 \times \FRAK{g}^2$
and
  $\FRAK{g}_1 = V_1 \otimes V_2$ where the $\FRAK{g}^i$ are
  nonzero semisimple Lie algebras and the $V_i$ are simple $\FRAK{g}^i$-modules.
  Furthermore, for $i = 1,2$ there are $\FRAK{g}^i$-invariant bilinear maps
  \[ \pi_i : V_i \times V_i \longrightarrow \FRAK{g}^i \quad ,
  \quad \psi_i : V_i \times V_i \longrightarrow \mathbb{C} \]
  such that
  \begin{equation} \label{e1}
   [u_1 \otimes u_2, v_1 \otimes v_2 ] =
  \psi_2(u_2,v_2)\pi_1(u_1,v_1) + \psi_1(u_1,v_1)\pi_2(u_2,v_2)
  \end{equation}
    for $u_1, v_1 \in V_1; u_2, v_2 \in V_2.$
In addition we can assume that $\pi_2, \psi_1$ are symmetric and
$\pi_1,
  \psi_2$ are skew-symmetric.

We claim that if $\FRAK{g} \neq \Gamma(\sigma_1, \sigma_2,
\sigma_3)$ there are nonzero constants $s_i$ such that
  the maps $\pi_1, \pi_2$ are given by
\begin{equation} \label{e70}
 \pi_i(u,v)w = s_i( \psi_i(v,w)u - \psi_i(w,u)v)
  \end{equation}
for $u_i, v_i \in V_i,$ cf. [Sch, page 144, equation (5.16)].\\

Indeed, from Section \ref{3new} equation (\ref{e70}) holds when
$\FRAK{g} = osp(m,n)$ with $m \geq 3, n \geq 2$.  Also  equation
(\ref{e70}) defines  $\FRAK{g}^i$-invariant bilinear maps $\pi_i :
V_i \times V_i \longrightarrow \FRAK{g}^i$ , so (\ref{e70})  holds
whenever $\FRAK{g}^1$ and $\FRAK{g}^2$ are simple and the adjoint
representation of $\FRAK{g}^1$
 resp. $\FRAK{g}^2$
occurs with multiplicity one in
  $\wedge^2V_1$ resp. $S^2V_2.$ This is the case for the Lie superalgebras $G(3)$ and $F(4).$
Note however that if
 $\FRAK{g} = \Gamma(\sigma_1, \sigma_2, \sigma_3)$
then we can write  $\FRAK{g}_0$ as  $\FRAK{g}^1 \times \FRAK{g}^2$
where $\FRAK{g}^1 \cong so(4)$ and  $\FRAK{g}^2 \cong sl(2).$ In
this case the  map $\pi_1 : V_1 \times V_1 \longrightarrow
\FRAK{g}^1$ is not, in general given by (\ref{e70}). This exhausts
all the classical simple Lie superalgebras  $\FRAK{g}$ such that
$\FRAK{g}_0$ is not simple and $\FRAK{g}_1$
  is an irreducible $\FRAK{g}_0$-module.\\

 Choose an orthonormal basis
  $\{ e_1, \ldots, e_m \}$ for
$V_1$ with
 respect to $\psi_1$ and a basis $\{ f_{\pm 1},
  \ldots, f_{\pm s}\}$ for $V_2$
 such that
  the matrix
 $J$ of $\psi_2$ with
 respect to this basis take the form
\[ J = \left[ \begin{array} {cc}
    0 & \Upsilon_s \\
 -\Upsilon_s & 0
   \end{array} \right] .\]

We denote by $so(V_1), sp(V_2)$ the orthogonal and symplectic
algebras preserving the forms $\psi_1, \psi_2$ respectively.  Let
$A_{ik} = \pi_1(e_i,e_k)$ and $B_{j \ell} = \pi_2(f_j,f_\ell)$.
 \begin{lemma}
 (a) With respect to the basis $\{e_i \otimes f_j\}$ of
$\FRAK{g}_1,$ we have
\[ M(\FRAK{g}) = A \otimes J + I_m \otimes B .\]
\\
(b)  Assume that the trace form $(a,b) \longrightarrow trace(ab)$
is used to evaluate matrices over $g\ell (V_1)$ and  $g\ell
(V_2).$ Then there are nonzero constants $\lambda, \mu$ such that
\[ A(x) = \lambda x, \; (JB)(y) = \mu y \]
for all $x \in so(V_1)$, $y \in sp(V_2) .$
\end{lemma}
 {\bf Proof.} (a) This
follows easily from formula $(\ref{e1})$.\\
\\
(b) We prove the statement about $so(V_1)$; the other part is
similar. We assume that $s_1 = 1.$ By equation (\ref{e70})
$\pi_1(e_i,e_j) = e_{i, j} - e_{j,i}.$ Hence \[ A = \sum_{i <
j}(e_{i, j} - e_{j,i})({\bf e_{i,j} - e_{j,i} }). \] The result
follows since for $i < j$ and $k < \ell$ we have
\begin{eqnarray*}
trace((e_{i, j} - e_{j,i})(e_{k, \ell} - e_{\ell,k})) =
-2\delta_{i,k}\delta_{j,\ell}.
 \end{eqnarray*}

\subsection {\label{3.4}}
To apply Lemma \ref{3.3} we need to consider three cases
separately. Suppose first that $\FRAK{g} = osp(m,n)$ with $m \geq
3, n \geq 2$. Then $\FRAK{g}_0 = \FRAK{g}^1 \times \FRAK{g}^2$
where $\FRAK{g}^1 = so(m), \FRAK{g}^2 = sp(n)$.  Also $\FRAK{g}_1
= V_1 \otimes V_2$ where $V_1$ is the natural module for $so(m)$
and $V_2$ is the natural
  module for $sp(n)$.  There are maps
  $\pi_i,  \psi_i$ for $i = 1,2$ such that the product
  $\FRAK{g}_1 \times \FRAK{g}_1 \longrightarrow \FRAK{g}_0$
  is given by equation (\ref{e1}) in section \ref{3.3}.

We recall how nilpotent orbits in simple Lie algebras of types
$B$, $C$ and $D$ can be described in terms of partitions. Let
${\bf P}_1(m)$ (resp. ${\bf P}_{-1}(m)$) be the set of partitions
of $m$ in which even (resp. odd) parts occur with even
multiplicity. Then by [CM, Theorems 5.1.2 and 5.1.3], nilpotent
orbits in $so(2r+1), r \geq 1$ (resp. $sp(2s), s \geq 1$) are in
one-to-one correspondence with partitions in ${\bf P}_1(2r+1)$
(resp. ${\bf P}_{-1}(2s))$.  We denote the orbit corresponding to
a partition $\mu$ by ${\cal O}_\mu$. We say that a partition is
{\em very even} if it has only even parts, each with even
multiplicity. By [CM, Theorem 5.1.4] any partition $\mu \in {\bf
P}_1(2r)$ corresponds to a unique orbit ${\cal O}_\mu$ in $so(2r),
r \geq 1$ unless $\mu$ is very even in which case $\mu$
corresponds to two orbits ${\cal O}^I_\mu$ and ${\cal
O}^{II}_\mu$.\\
 \indent
From the proofs of [CM, Propositions 5.2.3, 5.2.5 and 5.2.8],
 we see that if a simple Lie algebra of type $B$, $C$ or $D$ is
regarded as a subalgebra of $g \ell(N)$ using the defining
representation then the Jordan form of a matrix in ${\cal O}_\mu$
(or ${\cal O}^I_{\mu}, {\cal O}^{II}_{\mu}$) corresponds to the
partition $\mu$.

If $\mu \in {\bf  P}_1(m), \nu \in {\bf  P}_{-1}(2s)$ and $\mu$ is
not very even, we consider the orbits
\[ {\cal O}_{\mu,\nu} = {\cal O}_\mu \times {\cal O}_\nu .\]
If $m = 2r$ and $\mu$ is very even the existence of two orbits
${\cal O}^I_\mu$ and ${\cal O}^{II}_\mu$ causes some notational
difficulties.  The simplest solution is to abuse notation slightly
and allow ${\cal O}_{\mu,\nu}$ to denote either of the orbits
${\cal O}^I_\mu \times {\cal O}_\nu$ or ${\cal O}^{II}_\mu \times
{\cal O}_\nu$. Since the values of $k({\cal O}^I_\mu \times {\cal
O}_\nu)$ and $k({\cal O}^{II}_\mu \times {\cal O}_\nu)$ turn out
to be
the same this does not create any problems.
  \\
   \\
  \subsection {\label{3.5}}
Let $\FRAK{g} = G(3)$, then $\FRAK{g}_0 = \FRAK{g}^1 \times
\FRAK{g}^2$
  and $\FRAK{g}_1 = V_1 \otimes V_2$ where $\FRAK{g}^1 \cong  \FRAK{g}_2,$ the $14$
   dimensional
  exceptional simple Lie algebra,  $ \FRAK{g}^2 \cong
   s\ell(2),   V_1$ is the
  7-dimensional simple $\FRAK{g}_2$-module and $V_2$ the
  2-dimensional simple $s\ell(2)$-module.  There are invariant
  maps $\pi_1 : \wedge^2 V_1 \longrightarrow \FRAK{g}_2, \; \pi_2: S^2 V_2 \longrightarrow
  s \ell(2)$ and invariant bilinear forms $\psi_1, \psi_2$ such
  that the product $\FRAK{g}_1 \times \FRAK{g}_1 \longrightarrow
  \FRAK{g}_0$ is given by equation $( \ref{e1} )$.  In particular since
  $\FRAK{g}_2$ preserves $\psi_1$ it can be regarded as a
  subalgebra of $so(V_1) = so(7)$.  If $\cal O$ is a nilpotent orbit in
  $\FRAK{g}_2$ we write ${\cal O} = {\cal O}_\mu$ where $\mu$ is
  the partition of $7$ determined by the Jordan form of a
  representative element of ${\cal O}$ when viewed as an element
  of $g \ell (V_1)$.  These partitions, together with the usual
  Bala-Carter notation for orbits in $\FRAK{g}_2$ [CM, page 128] and the dimension of
  the orbits are given in the table below.\\
  \\
  \[  \begin{tabular}{|c|c|c|c|c|c|} \hline
 ${\cal O}={\cal O}_\mu$&${0}$ & $A_1$ & $\widetilde{A}_1$
             & $G_2(a_1)$ & $G_2$ \\ \hline
 $\mu$ & $1^7$ & $2^2,1^3$ & $3,2^2$ & $3^2, 1$ & $7$  \\ \hline
 $\dim{\cal O}$ & $0$ & 6 & 8 & 10 & 12\\ \hline
\end{tabular} \]
For $\mu$ in the table and $\nu \in {\bf  P}(2)$ set  ${\cal
O}_{\mu,\nu} = {\cal O}_\mu \times {\cal O}_\nu .$
  \\

  In section \ref{3.7} we apply Lemma \ref{3.3} to calculate $k({\cal
  O}_{\mu,\nu}).$  However to do this we need to evaluate the
  matrix using an invariant bilinear form on $\FRAK{g}_2,$
rather than on $g\ell (7)$
  as was done in Lemma \ref{3.3}.  Similar remarks apply when
  $\FRAK{g}$ is the  Lie superalgebra F(4).  Recall that any nonzero invariant form on a simple Lie algebra is proportional
  to the Killing form.  Therefore since  ${\cal
  O}_{\mu,\nu}$ is a product of conical subvarieties (see section
  \ref{EM}), our method is justified by the following well-known
  lemma.  Our proof is a modification of [LS, Lemma 2.5].
  \begin{lemma} {\label{jl}}
Suppose that $\mathfrak k \subseteq \mathfrak l$ are finite
dimensional simple complex Lie algebras.  Then the restriction of
the Killing form $B$ on $\mathfrak l$ to $\mathfrak k$ is
nondegenerate.
\end{lemma}
{\bf Proof.} There are connected, simply connected complex Lie
groups $K$ and $L$, unique up to isomorphism, such that
$\mathfrak k = Lie(K)$ and $\mathfrak l = Lie(L)$. We can take $K$
to be a subgroup of $L$ since $\mathfrak k \subseteq \mathfrak l$.

Let $K_0$ denote a maximal compact subgroup of $K$. Then  $K_0$ is
contained in a maximal compact subgroup, $L_0$, of $L$.  Let
$\mathfrak k_0$ (resp.  $\mathfrak l_0$) denote the (real) Lie
algebra of the compact Lie group, $K_0$ (resp. $L_0$).  We have
$\mathfrak l = \mathfrak l_0 \oplus i \mathfrak l_0$ and
$\mathfrak k = \mathfrak k_0 \oplus i \mathfrak k_0$.

Now $B$ is negative definite when restricted to $\mathfrak l_0$
and hence it is negative definite on $\mathfrak k_0.$ Therefore
the restriction $B'$ of $B$ to $\mathfrak k$ is nonzero.  However
the radical of $B'$ is an ideal of $\mathfrak k,$ so $B'$ is
nondegenerate.
  \\
\subsection {\label{3.6}}
 Now let $\FRAK{g} = F(4)$.  Then  $\FRAK{g}_0
= \FRAK{g}^1 \times \FRAK{g}^2$
  and $\FRAK{g}_1 = V_1 \otimes V_2$ where $\FRAK{g}^1 \cong  so(7),
   \FRAK{g}^2 \cong
   s\ell(2),  V_1$ is the spin
representation of $so(7)$ and $V_2$ is the 2-dimensional simple
$s\ell(2)$-module. We have the same analysis as for $G(3)$ except
that $so(7)$ is now regarded as a subalgebra of $so(V_1) = so(8)$.

Nilpotent orbits in $so(7)$ correspond to partitions $\eta \in
{\bf  P}_1(7)$.  For $\eta \in {\bf  P}_1(7)$
  we write $\mu = \sigma(\eta)$ for where $\mu$ is
  the partition of $8$  determined by the Jordan form of an
element of the corresponding orbit when viewed as an element
  of $g \ell (V_1)$.  We use $\mu$ to label the orbit.
  The map $\sigma:{\bf P}_1(7) \longrightarrow {\bf P}(8)$, together with the
  dimension of the orbits are given in the table below.
  \\
  \\
  \[  \begin{tabular}{|c|c|c|c|c|c|c|c|} \hline
 $\eta$ & $1^7$&$2^2,1^3$ & $3,1^4$ & $3,2^2$
             & $3^2,1$ & $5,1^2$ & $7$ \\ \hline
 $\mu = \sigma(\eta)$ & $1^8$ & $2^2,1^4$ & $2^4$ & $3,2^2,1$ & $3^2,1^2$ & $4^2$ & $7,1$\\ \hline
 $\dim{\cal O_\mu}$ & $0$ & $8$ & $10$ & $12$& $14$& $16$ & $18$\\ \hline
\end{tabular} \]

As before we set  ${\cal O}_{\mu,\nu} = {\cal O}_\mu \times {\cal
O}_\nu $ for $\mu$ in the table and $\nu \in {\bf  P}(2).$\\
 \\

\subsection {\label{3.7}}
  Let $\FRAK{g} = osp(m,n)\; (m \geq 3),\;  G(3)$ or $F(4)$ and consider the nilpotent orbit
 ${\cal O}_{\mu,\nu}$ as defined in one of the three preceding subsections.
Let $\dim V_1  = m$ and $\dim V_2 = n.$

  \begin{theorem}
We have
\[ k({\cal O}_{\mu,\nu}) = dim \FRAK{g}_1  - \sum_{i}\mu'_i \nu'_i. \]
\end{theorem}
{\bf Proof.} If $(x,y) \in {\cal O}_{\mu,\nu} $ then $k(q)$ is the
rank of the evaluation of $M(\FRAK{g})$ at $(x,y)$.  Since  $I_m
\otimes J$ is invertible this is the same as the rank of the
evaluation of $(I_m \otimes J)M(\FRAK{g}) = -A \otimes I_n + I_m
\otimes B$ at $(x,y)$.  Thus the result follows from Lemmas
\ref{Part} and \ref{3.3} .\\

 \subsection {\label{3.8}} Theorem \ref{3.7} does not apply to the Lie superalgebras
$\FRAK{g} = osp(m,2r)$ when $m = 1, 2$. To handle these cases we
use the description of $osp(m,n)$ given in section \ref{3new}.

If $m = 1,$ we choose $e \in V_1$ such that $\psi_1(e,e) = 1.$
Then for $v, w \in V_2$ we have
 \begin{equation} \label{e60}
 [e \otimes v, e \otimes w] = v \circ w.
  \end{equation}

If $m = 2,$ we choose $e_{-}, e_{+} \in V_1$ such that
\[ \psi_1(e_{-}, e_{-}) = \psi_1(e_{+}, e_{+}) = 0,  \quad \psi_1(e_{-}, e_{+}) = 1. \]
Set   $ \FRAK{g}_1^{\pm} = \mathbb{C}e_{\pm} \otimes V_2,$ and $z
=  e_{-} \wedge e_{+}.$ Then $\FRAK{g}_0 = [\FRAK{g}_0,
\FRAK{g}_0] \oplus \mathbb{C}z,$ and
 $\FRAK{g}_1 = \FRAK{g}_1^+ \oplus \FRAK{g}_1^-,$
is a direct sum of $\FRAK{g}_0$-modules. Also  $
[\FRAK{g}_1^{\pm},\FRAK{g}_1^{\pm}] = 0$ and for $v, w \in V_2$ we
have
 \begin{equation} \label{e50}
 [e_{-} \otimes v, e_{+} \otimes w] = v \circ w + \psi_2(v,w)z.
  \end{equation}

 If  $\FRAK{g} = osp(m,2r)$ where $m = 1,
2$, then nilpotent orbits in $ \FRAK{g}_0$ are parameterized by
partitions in
${\bf P} _{-1}(2r)$. We denote the orbit corresponding to a
partition $\mu$ by ${\cal O}_\mu$. Note that the rank of $J_{\mu}$
is $\sum_{i}(\mu_i -1) = 2r - \mu'_1.$

 \begin{theorem}

(a) If $\FRAK{g} = osp(1,2r)$ and $\mu \in {\bf  P}_{-1}(2r)$ we
have
\[ k({\cal O}_{\mu}) =  rank \; J_{\mu}. \]
(b) If $\FRAK{g} = osp(2,2r)$ and $\mu \in {\bf  P}_{-1}(2r)$ we
have
 \[ k({\cal O}_{\mu}) = 2(rank \; J_{\mu}). \]
 \end{theorem}
 {\bf Proof.} (a) Identify $ \FRAK{g}_1 =  \mathbb{C}e \otimes
 V_2$ with $V_2$ via the map $ e \otimes v \longrightarrow v.$
Let $K = \{ \pm 1, \ldots, \pm r \}$ and choose a basis $\{ e_k |
k \in K \}$ for $ \FRAK{g}_1$
  such that the matrix of $\psi$ on this basis is the matrix
$J$ used in the proof of Lemma \ref{3.3}.

The matrix $M(\FRAK{g})$ equals   $\sum_{i,j \in K} ( e_i \circ
e_j) {\bf e_{i,j}}  ,$ and as in the proof of Lemma \ref{3.3}
there is a nonzero constant $\lambda$ such that
\[ JM(\FRAK{g})(x)  = \lambda x\]
for all $x \in \FRAK{g}_0.$  This easily gives the result.\\
 (b) Let $\FRAK{g} = osp(2,2r),$ and $\FRAK{k} = osp(1,2r).$
By comparing equations (\ref{e60}) and (\ref{e50}), we see that
with respect to a suitable ordered basis, $M(\FRAK{g})$ has the
block matrix form
\[  \left[ \begin{array}{cc}
    0     & M(\FRAK{k})\\
   M(\FRAK{k}) & 0
  \end{array} \right]  \quad \mbox{mod} \; (z). \]
The result follows since $z$ vanishes on any nilpotent orbit in
$\FRAK{g}_0$.

\subsection {\label{3.9}}
Now let $\FRAK{g} = \Gamma(\sigma_1, \sigma_2, \sigma_3)$ as in
[Sch]. Then $\FRAK{g}_0 = \FRAK{g}^1 \times \FRAK{g}^2 \times
\FRAK{g}^3,\; \FRAK{g}_1 = V_1 \otimes V_2 \otimes V_3$ where
$\FRAK{g}^i \cong s\ell(2) $ and $V_i$ is the 2-dimensional simple
$s\ell(2)$-module.

 Let $\psi_i : V_i \times V_i \longrightarrow  \mathbb{C}$
be a nonzero
 $\FRAK{g}^i$-invariant skew-symmetric map and
define a $\FRAK{g}^i$-invariant symmetric map

\[ \pi_i: V_i \times V_i \longrightarrow
 \FRAK{g}^i \]
by
\[ \pi_i(x,y)z = \psi_i(y,z)x - \psi_i(z,x)y \]
for $x,y,z \in V_i$.
  Then for $a_1 \otimes a_2 \otimes a_3, b_1 \otimes b_2 \otimes b_3 \in \FRAK{g}_1$
   we have

\begin{equation}
\label{G2}
 [a_1 \otimes a_2 \otimes a_3, b_1 \otimes b_2 \otimes b_3] = \\
= \sum \sigma_k \psi_i (a_1, b_1)\psi_j (a_2,b_2)\pi_k (a_3,b_3)
  \end{equation} \\
where the sum is over all even permutations $(i,j,k)$ of
$\{1,2,3\}$.  Let $f,h,e$ be the basis of $s \ell(2)$ given by
\[ f = \left[ \begin{array}{cc}
 0 & 0 \\
 1 & 0 \end{array}
 \right], \quad
 h = \left[ \begin{array}{cc}
1 & 0 \\
 0 & -1 \end{array}\right], \quad
 e = \left[ \begin{array}{cc}
  0 & 1 \\
  0 & 0 \end{array}\right] \]
 and let $x = (1, 0)^t$ and $y = (0, 1)^t$
  be basis vectors for the
  2-dimensional $s\ell(2)$-module.  We write $f_i, h_i, e_i$
  (resp. $x_i, y_i$) for the corresponding elements of
  $\FRAK{g}^i$ (resp. $V_i)$, and set $S_i = \{x_i, y_i\}$.
Consider the matrices
\[ \Psi = \left[ \begin{array}{cc}
 0 & 1 \\
 -1 & 0 \end{array}\right], \quad
 \Pi_i = \left[ \begin{array}{cc}
2e_i & -h_i \\
 -h_i & -2f_i \end{array}\right]. \]

We assume  that the matrix for each $\psi_i$ on the ordered
   basis $(x_i, y_i)$ for $V_i$ is
$\Psi.$  Then the matrix for $\pi_i$ on this
   basis is $\Pi_i.$
   We order the basis $\{a_1
   \otimes a_2 \otimes a_3|a_i \in S_i \}$ of $\FRAK{g}_1$
   lexicographically.  It follows from equation (\ref{G2}) that the matrix $M(\FRAK{g})$
is given by
\begin{equation}
\label{G3}
 M(\FRAK{g}) = \sigma_3 \Psi \otimes \Psi \otimes \Pi_3  + \sigma_2
   \Psi \otimes \Pi_2 \otimes \Psi  + \sigma_1 \Pi_1 \otimes \Psi \otimes \Psi
  \end{equation} \\
This can also be deduced from  Table I in [Z].

For $ \mu, \upsilon ,\eta \in {\bf P}(2)$ let $\cal O_{\{ \mu,
\upsilon ,\eta \} }$ denote the orbit of    $(J_{\mu}, J_{
\upsilon} ,J_{\eta })$ in $\FRAK{g}_0.$ Note that the evaluation
of the matrix $\Psi^{-1} \Pi_i $ at any element $x$ of
$\FRAK{g}^i$ is a nonzero multiple of $x.$ It follows from
equation (\ref{G3}) that we can find $x \in \cal O_{\{ \mu,
\upsilon ,\eta \} }$ such that the evaluation of $(\Psi \otimes
\Psi \otimes \Psi)^{-1} M(\FRAK{g})$ at $x$ equals
\[ J_\mu \otimes I_2  \otimes I_2 + I_2 \otimes J_\nu  \otimes I_2 + I_2  \otimes  I_2 \otimes
J_\eta .\] The values of $\dim{\cal O_{\{ \mu, \upsilon ,\eta
\}}}$ and
 $k({\cal O_{\{ \mu, \upsilon ,\eta \}}})$ depend only on the set $\{ \mu, \upsilon ,\eta \}.$
These values are given in the table below.
 \[  \begin{tabular}{|c|c|c|c|c|} \hline
 $\{ \mu, \upsilon ,\eta \}$&$\{ 2,2,2 \}$ & $\{ 2,2,1^2 \}$ & $\{ 2,1^2,1^2 \}$
             & $\{ 1^2,1^2,1^2 \}$ \\ \hline
 $\dim{\cal O_{\{ \mu, \upsilon ,\eta \}}}$ & $6$ & $4$ & $2$ & $0$  \\ \hline
 $k({\cal O_{\{ \mu, \upsilon ,\eta \}}})$ & $5$ & 4 & 4 & 0 \\ \hline
\end{tabular} \]
We can view $\FRAK{g}$ as a deformation of $D(2,1) = osp(4,2)$ and
the values of $k(\cal O)$ for $\FRAK{g}$ are the
same as those for the corresponding orbits for $D(2,1)$.\\

\subsection {\label{3.11}}
Let $V_0,V_1$ be vector spaces with bases $e_1, \ldots, e_n$ and
$e'_1, \ldots e'_n$ respectively, and let $\psi : V_0
\longrightarrow V_1$ be the map sending $e_i$ to $e'_i$ and $e'_i$
to $-e_i$.  Let $\FRAK{g}$ denote the Lie superalgebra of all
endomorphisms of $V = V_0 \oplus V_1$ which supercommute with
$\psi$. Then $\FRAK{g}$ is isomorphic to the Lie superalgebra of
matrices of the form
\[ \left[ \begin{array}{cc}
    a & b \\
    b & a
    \end{array}\right] \]
    with $a,b \in g \ell (n)$.  Thus $\FRAK{g}_0 \cong g \ell(n)$
    and $\FRAK{g}_1 \cong \FRAK{g}_0$ as a $\FRAK{g}_0$-module.
    The derived algebra $\FRAK{g}'$ consists of all matrices as
    above with $b \in s \ell(n)$. Also $\FRAK{g}'$ has a
    one-dimensional center $\FRAK{z} = \mathbb{C}I_{2n}.$  The factor algebra $\FRAK{g}'/\FRAK{z}$ is
    the simple Lie superalgebra denoted $Q(n-1)$ in [K].  We assume
    that $n \geq 3$.  Then the Lie superalgebra $Q(n-1)$ is simple.
    As a Cartan subalgebra $\FRAK{h}_0$ of $\FRAK{g}_0$ we take all matrices of the above
    form with $a$ diagonal and $b = 0$.  We modify this in the obvious
    way to obtain Cartan subalgebras of $\FRAK{g}'_0$ and $Q(n-1)_0$.

    \indent If ${\cal O}$ is any nilpotent orbit in $\FRAK{g}'_0$
    then $\FRAK{z}$ vanishes on ${\cal O}$ and ${\cal O}$ may be
    regarded as a nilpotent orbit in $(\FRAK{g}'/\FRAK{z})_0$.  All
    nilpotent orbits in $(\FRAK{g}'/\FRAK{z})_0$ arise in this
    way.  Therefore it suffices to consider the Clifford algebras
    arising from $\FRAK{g}$ and $\FRAK{g}'$.

    \indent If $\mu \in {\bf  P}(n)$ let $J_\mu$ and
     ${\cal O}_\mu$ denote the corresponding Jordan matrix and
     nilpotent orbit.  Set
      \begin{eqnarray*}
      \epsilon(\mu) & = & 1 \quad \mbox{if all parts of}\; \mu \; \mbox{are even} \\
                    &   & 0 \quad \mbox{otherwise}.
       \end{eqnarray*}
Note that $\FRAK{g}'_0 = \FRAK{g}_0$.  If
 $q \in Spec S(\FRAK{g}_0)$ let $k(q)$ (resp. $k'(q))$ be the rank
 of the bilinear form on $\FRAK{g}_1 \otimes F_q$ (resp
 $\FRAK{g}^{'}_{1} \otimes F_q)$ defined in the usual way.\\
 \\
\begin{theorem} \label{Q0}
 If $V(q) = \overline{\cal O}_\mu$ then \\
 \\
 (a) $k(q) = \dim \FRAK{g}_1 - \sum_i(\mu'_i)^2$\\
 \\
 (b) $k'(q) = k(q) - 2 \epsilon (\mu).$\\
\end{theorem}
{\bf Proof.} For $a \in g\ell(n)$ set
 \[ \overline{a} = \left[
 \begin{array}{cc}
 0               & a \\
 a  & 0
   \end{array}\right] . \]
 Let $K = \{ 1, \ldots, n \}$ and calculate $M(\FRAK{g})$
 using the basis $\{ \overline{e}_{ij} \}$ of $\FRAK{g}_1$.  The
 rows and columns of $M(\FRAK{g})$ are indexed by $K \times
 K$ ordered lexicographically with entry in row $(i,j)$ and column $(k,\ell)$ given by
 \[ [\overline{e}_{ij}, \overline{e}_{k\ell}] = \delta_{jk}
 e_{i\ell} + \delta_{i \ell} e_{kj} . \]
Thus
\[ M(\FRAK{g}) = \sum_{i,j,k,\ell}(\delta_{jk}e_{i\ell} +
\delta_{i\ell}e_{kj})  {\bf e_{i,k}} \otimes {\bf e_{j,\ell } }
.\] If $L = \sum {\bf e_{r,s}} \otimes {\bf e_{s,r}},$  then $L$
is nonsingular since   $L^2$ is the identity matrix. Let $A =
\sum_{i,j \in K} e_{ij} {\bf e_{i,j}}.$
 Then  $({\bf e_{i,k}} \otimes {\bf e_{j \ell}}) L = {\bf
 e_{i,\ell}}
\otimes {\bf e_{j,k}}$
 and hence
 \[ M(\FRAK{g}) L = A \otimes I_n + I_n \otimes A^{t} .\]
\\
Since $M(\FRAK{g})$ and $M(\FRAK{g}) L$  have the same rank, part
(a) of the Theorem follows from Lemma \ref{Part}.
Part (b) follows from the lemma in the next subsection.\\

\subsection {\label{3.12}}
With $\FRAK{g},\FRAK{g}'$ as in subsection \ref{3.11} we compare
the matrices $M(\FRAK{g})$ and $M(\FRAK{g}')$. For $1 \leq i \leq
n-1$ let $h_i = e_{ii} - e_{i+1, i+1}$ and let $h_n$ be the
identity matrix. We calculate $M(\FRAK{g})$ using the basis
 \[ \{ \overline{e}_{ij}, \overline{h}_k | 1 \leq i \neq j \leq n, 1 \leq k \leq n \} \]
  of
 $\FRAK{g}_1.$  We order this basis in any way such that the last
 $n$ elements are $\overline{h}_1, \ldots , \overline{h}_n .$

 Note that for $1 \leq i \leq n - 1$ we have
\begin{equation} \label{e4}
[ \overline{h}_i, \overline{e}_{k,k+1}] = (\delta_{i,k+1} - \delta_{i,k-1} )e_{k,k+1}\\
\end{equation}

\begin{equation}
\label{e5}
 [ \overline{h_n}, \overline{e}_{k,k+1}] = 2e_{k,k+1}.
  \end{equation}

  The evaluation of $M(\FRAK{g})$ at $J_\mu$ has the block-matrix
  form
  \[ \left[ \begin{array}{c|c}
   * & N(J_{\mu}) \\ \hline
    N(J_{\mu})^{t} & 0
    \end{array}\right] \]
    where $N$ is the matrix with entries $[\overline{h}_i,
    \overline{e}_{k \ell}]$ $( 1 \leq i \leq n, 1 \leq k \neq \ell
    \leq n)$.  The evaluation of $M(\FRAK{g}')$ at $J_\mu$ is
    obtained by deleting the last row  and column.

    \indent For $i \neq j$ let $C_{ij}$ be the column of
    $M(\FRAK{g})$ corresponding to $\overline{e}_{ij}$.  Also
    for $1 \leq i \leq n-1$ let $C_i$ be the column corresponding
    to $\overline{h}_i$.  The evaluation of a column $C$ at
    $J_\mu$ is denoted $C(J_\mu)$.
      \\
      \\

     \begin{lemma}
(a) If $(\sum_{k,\ell} \lambda_{k,\ell} C_{k,\ell} + \sum^n_{k=1}
\nu_k
    C_k)(J_\mu) = 0$
     then
    \[ \sum^n_{i=1} \nu_i C_i (J_\mu) = 0 .\]

(b) The linear span of the columns $C_1(J_{\mu}), \ldots,
C_{n-1}(J_\mu)$ contains $C_n(J_\mu)$ if and only if some part of
$\mu$ is odd.
\end{lemma}
{\bf Proof.} (a) This follows since $(\sum_k \nu_k C_k)(J_\mu)$
can have nonzero entries only in rows $(i, i+1)$ and $(\sum
\lambda_{k\ell} C_{k\ell})(J_\mu)$ has zero entries in these
rows.\\
 (b) Consider the system of equations
  \begin{equation}\label{e8}
    C_n(J_\mu) = 2 \sum^{n-1}_{i=1} x_i C_i(J_\mu)
 \end{equation}
 in the unknowns $x_1, \ldots, x_{n-1}.$
  By equations (\ref{e4}) and (\ref{e5}) this system is equivalent to the
  evaluation of the system of equations
\begin{equation} \label{e6}
 e_{k,k+1} = (x_{k+1}  - x_{k-1}) e_{k, k+1}
  \end{equation}
  at $J_\mu$.  Here we set $x_0 = x_n = 0$.  Thus the system
  (\ref{e8}) is equivalent to the equations
   \begin{equation}\label{e9}
   1 = x_{k+1} - x_{k-1} \; \mbox{for} \; 1 \leq k \leq n-1, k \neq \mu_1 +
   \ldots + \mu_i.
   \end{equation}
   If $\mu_i$ is even for all $i,$ then (\ref{e9}) involves the
   equations
   \[ 1 = x_n - x_{n-2} = \ldots = x_2 - x_0 \]
   which are inconsistent.

On the other hand if some $\mu_i$ is odd, then $\mu_1 +\ldots +
\mu_j$ is odd for some $j$, so the system (\ref{e9}) is equivalent
to a number of systems of equations of the form
\begin{equation} \label{e7}
1 = x_{p} - x_{p-2} = \ldots = x_{q+2} - x_{q}.
\end{equation}
Moreover the sets of variables which occur in two such systems are
disjoint, and in each system (\ref{e7}) we have either $p < n$ or
$q
> 0$. If $p < n$ (resp. $q
> 0$), we can set $x_q = 0$ (resp. $x_p = 0$)  and solve the
equations (\ref{e7}) recursively for $x_{q+2i}$ (resp.
$x_{p-2i}$).
\\
\\
{\bf Remark.}  If $V(q) = \overline{\cal{O}}_\mu$ it follows from
Theorem \ref{3.11} and [CM, Corollary 7.2.4] that $k(q) = \dim
\overline{\cal{O}}_\mu$.

\subsection {\label{3.10}}
For any classical simple Lie superalgebra $\FRAK{g}$ considered up
to this point, the matrix $M(\FRAK{g})$ is nonsingular.  This fact
together with some Clifford algebra theory can be used to show
that $U(\FRAK{g})$ is prime, [B].  However if $\FRAK{g} = P(n),$
it is shown in [KK] that $U(\FRAK{g})$ is not prime, and it
follows that $M(\FRAK{g})$ is singular.  Because of this it seems
unlikely that $M(\FRAK{g})$ can be expressed in terms of a
Kronecker product.  However if $\mathcal{O}$ is a nilpotent orbit
in $\FRAK{g}_0$ there is a formula for $k(\mathcal{O})$ which is
similar to the formula for the corresponding orbit for the Lie
superalgebra $Q(n).$\\
\indent
For $n \geq 2$ the Lie superalgebra $P(n)$ is the
subalgebra of $s \ell(n + 1, n + 1)$  consisting of all matrices
of the form
\[ \left[ \begin{array}{cc}
A & B \\ C & -A^t
\end{array} \right] \]
where $\mbox{trace}\,(A) = 0, \, B^t = B$ and $C^t = - C$.

If $\FRAK{g} = P(n)$, then $\FRAK{g}_0 \cong s \ell(n+1).$ As a
$\FRAK{g}_0$-module, $\FRAK{g}_1$ is the direct sum of two
submodules $\FRAK{g}^{\pm}_1$ where $\FRAK{g}^+_1$ (resp.
$\FRAK{g}^-_1$) consists of all matrices as above with $B = 0$
(resp. $C = 0$).  Let $V$ be the natural module for $ s \ell(n+1)$
with weights  $\epsilon_1, \ldots ,\epsilon_{n+1}.$
 Then, as $\FRAK{g}_0$-modules
 $\FRAK{g}_1^+ \cong S^2 V$ and
$ \FRAK{g}_1^{-} \cong \wedge^2 V^*.$\\

 Fix $\mu$ a partition of $n$.  We
assume that the nonzero entries in the Jordan matrix $J_\mu$ occur
immediately below the main diagonal. For $1 \leq i \leq n - 1$,
let $b_i$ be the entry of $J_\mu$ in row $i+1$, column $i$, and
let  $b_0 = b_n = 0$.  Denote the orbit of ${\cal O}_{\mu}$ in
$\FRAK{g}_0$ by $J_\mu.$

\begin{theorem}
For $\mu \in {\bf  P}(n)$  we have
\[ k({\cal O}_{\mu}) = 2 \sum^{n-1}_{i=1}(n-i)b_i = n^2 - \sum_i(\mu'_i)^2 . \]
\end{theorem}
 {\bf Proof.}  If we choose a basis for
 $\FRAK{g}_1$ such that elements of
$\FRAK{g}_1^+$ precede elements of $\FRAK{g}_1^{-},$ then
$M(\FRAK{g})$
 has the form
 \[ \left[ \begin{array}{cc}
 0 & N \\
 N^{t} & 0
 \end{array}\right].\]
Let $\epsilon_1, \ldots ,\epsilon_{n+1}$ be the weights of
 $V$.  We use the weights $-\epsilon_i - \epsilon_j$ of
 $\wedge^2 V^*$ to index the rows, and the weights $\epsilon_k +
 \epsilon_\ell$ $(k \leq \ell)$ of $S^2 V$ to index the columns of $N$.
 It is easy to see that if
 if $x \in \wedge^2 V^*$ has weight $- \epsilon_i - \epsilon_j$
 and $y \in S^2 V$ has weight $\epsilon_k + \epsilon_\ell$  and $(i,j) \neq (k, \ell )$ then
 $[x,y] \neq 0.$  Hence
 the entry in row $(i,j)$ and column $(k,\ell)$ of $N$
 is, up to a nonzero scalar equal to $[x,y]$.
We order the rows of $N$ lexicographically and  order the  columns
so that column $(i,j)$  precedes  column $(k,\ell)\; \mbox{if and
only if} \; (k, \ell) <_{\ell ex}(i,j).$

\indent Note that the evaluation $N_\mu$ of $N$ at
$J_\mu$ has the following properties:\\
\\
(a) The entry in row $(i,j)$ and column $(i + 1, j)$ is nonzero if and only if  $b_i = 1$ \\
\\
(b) The entry in row $(i,j)$ and column $(i, j+1)$ is nonzero if and only if $b_j = 1$.\\
\\
(c) All other entries in row $(i,j)$ are zero.
\\
\\
  We claim that $N_\mu$ is row equivalent to the matrix
$\overline{N}_\mu$ obtained from $N_\mu$ by replacing row $(i,j)$
by zero for all $j > i$ whenever $b_i = 0$. We can assume that $i
> 1$, since if $b_1 = 0$ then $J_\mu = 0$, and also that $b_j = 1$.
Then $b_{i-1} = 1$.  Suppose that  $b_i = b_{i-q-1} = 0$ but
$b_{i-p} \neq 0$ for $p = 1, \ldots , q.$  This means that the
Jordan block of $J_\mu$ ending in row $i$ has size $q+1$.  Since
the Jordan blocks of $J_\mu$ are arranged in order of decreasing
size and $i < j$ it follows that the Jordan block of $J_\mu$
containing row $j$ has size at most $q + 1$.  Hence $b_{j+p} = 0$
for some $p$ with $p \leq q$ and we fix $p$ minimal with this
property. \indent Then the submatrix of $N_\mu$ formed by rows $(i
- p - 1 + k, j + p + 1 - k)$, for $k = 0, \ldots, p + 1$, and
columns $(i - k, j + k + 1),$ for $k = 0, \ldots, p$ has the form
 \[ \left[
\begin{array}{lll}
    0\;0\;0           & \cdot\cdot\cdot  & 0\;0\;\ast \\
    0\;0\;0          & \cdot\cdot\cdot  & 0\ast\ast\\
    0\;0\;0           & \cdot\cdot\cdot  & \ast\ast0 \\
  \cdot\cdot\cdot &                  & \cdot\cdot\cdot\\
  \cdot\cdot\cdot &                  & \cdot\cdot\cdot\\
   0 \ast \ast    & \cdot\cdot\cdot  & 0\;0\;0 \\
 \ast \ast 0      & \cdot\cdot\cdot  & 0\;0\;0 \\
    \ast \;0\; 0      & \cdot\cdot\cdot  & 0\;0\;0  \end{array} \right] \]
 where each * is nonzero.
In addition every nonzero entry in each of the rows listed above
occurs in this submatrix.  Hence the last row, row $(i,j)$ of
$N_\mu$, is a linear combination of the preceding rows.  The claim
follows from this.

Now if $b_i \neq 0$, then for $j > i$ the first entry in row
$(i,j)$ of $\overline{N}_\mu$ occurs in column $(i+1, j)$.  Each
such index $i$ contributes $n  - i$ linearly independent rows to
the rank of $\overline{N}_\mu$, so we obtain the first formula in the theorem.\\
\\
To obtain the second formula, note that $b_i = 0$ if $i = \mu_1 +
    \ldots + \mu_k$ for some $k$ and that $b_i = 1$ otherwise.
    Hence
    \begin{eqnarray*}
     k({\cal O}_\mu) &=& 2 \sum_i(n - i)b_i \\
    & = & n(n-1) - 2 \sum_k (n - (\mu_1 + \ldots + \mu_k)).
    \end{eqnarray*}
    Observe that $n - (\mu_1 + \ldots + \mu_k)$ is the number of
    boxes in the Young diagram for $\mu$ which are not contained
    in the first $k$ rows.  Using the columns instead to count boxes we have
    \begin{eqnarray*}
     k({\cal O}_\mu) & = & n(n-1) - 2 \sum_{i \geq 1} \sum_{j \geq
     1} (\mu'_i - j) \\
     & = & n(n-1) - \sum_{i \geq 1} \mu'_i(\mu'_i - 1) \\
      & = & n^2 - \sum_i(\mu'_i)^2.
    \end{eqnarray*}

\section{Parabolic Subalgebras}
\subsection {\label{4.1}}
\indent Although the connection with Clifford algebras works best
for the Lie superalgebras $g\ell(m,n)$, and $Q(n)$ many of our
results on induced modules hold more generally.  Therefore we
adopt an axiomatic approach.  Henceforth we assume that\\
\\
(i)  $\FRAK{g} = \oplus_{i \in \mathbb{Z}} \FRAK{g}(i)$ is a
graded
Lie superalgebra with $\FRAK{g}_0$ reductive.\\
\\
(ii)  $\FRAK{h}_0 \subseteq \FRAK{g}(0)$ where $\FRAK{h}_0$ is a
Cartan subalgebra (CSA) of $\FRAK{g}_0$ and $\FRAK{g}$ is a
semisimple $\FRAK{h}_0$-module.\\
\\
 Assume axioms (i) - (ii) and set
 \[ \FRAK{m} = \oplus_{i < 0} \FRAK{g}(i), \quad  \FRAK{l} =
 \FRAK{g}(0),\quad \FRAK{m}^+ = \oplus_{i > 0} \FRAK{g} (i), \quad \FRAK{p} =  \FRAK{l} \oplus
 \FRAK{m}^+  \]
 so that
 \[ \FRAK{g} = \FRAK{m} \oplus  \FRAK{l} \oplus
 \FRAK{m}^+ .\]

Let  $\FRAK{h}$ be the centralizer of $\FRAK{h}_0$ in $ \FRAK{g}.$
Axiom (ii) implies the existence of a root space decomposition
\[ \FRAK{g} = \FRAK{h} \oplus \oplus_{\alpha \in \Delta}
\FRAK{g}^\alpha \] where
\[ \FRAK{g}^\alpha = \{x \in \FRAK{g} | [h,x] = \alpha (h)x \;\;
\mbox{for all} \;\; h \in \FRAK{h}_0 \} \]
 and
 \[ \Delta = \{ \alpha \in \FRAK{h}_0^*|  \alpha \neq 0, \FRAK{g}^\alpha \neq 0 \} .\]
\indent

We also assume that\\
 (iii) $\Delta = \Delta^+ \cup \Delta^-$, a disjoint union, where $\Delta^\pm$ are
 subsets of $\Delta$ such that $\alpha, \beta \in \Delta^{\pm}$
 implies that $\alpha + \beta \in \Delta^\pm$ or $\FRAK{g}^{\alpha
 + \beta} = 0$, and  such that $\FRAK{g}^{\alpha} \cap  \FRAK{g}_0 \subseteq
\FRAK{p}$ for all $\alpha \in \Delta^+.$\\

Now let $\Delta(\FRAK{l})$ be the set of roots of $\FRAK{l}$ and
set $\Delta^{\pm}(\FRAK{l}) = \Delta^{\pm} \cap \Delta(\FRAK{l})
.$  If $\Gamma$ is a subset of $\Delta$ and $i = 0, 1$ we set
$\Gamma_i = \{\alpha \in
\Gamma| \FRAK{g}^{\alpha} \cap  \FRAK{g}_i \neq 0 \}.$\\
We refer to the subalgebra
\[ \FRAK{b} = \FRAK{h} \oplus \oplus_{\alpha \in \Delta^+}
\FRAK{g}^\alpha \]
 as a {\it Borel subalgebra} of $\FRAK{g}$.  Note that $\FRAK{b}$ is determined by $\Delta^+$ in axiom (iii)
  and
that in general there may be several choices for $\Delta^+$ even
if $\Delta_0^+$ is specified in advance.  By axiom (iii)
$\FRAK{b}_0 \subseteq \FRAK{p}.$   The subalgebra
\[ \FRAK{c} = \FRAK{h} \oplus \oplus_{\alpha \in \Delta^+(\FRAK{l})} \FRAK{g}^\alpha
= \FRAK{b} \cap  \FRAK{l} \] is a Borel subalgebra of $\FRAK{l}.$
We say that a root $\alpha \in \Delta^+(\FRAK{l})_0$ (resp.
$\alpha \in \Delta^+(\FRAK{l})$) is {\it indecomposable} if we
cannot write $\alpha$ in the form $\alpha' + \alpha''$ with
$\alpha', \alpha'' \in \Delta^+(\FRAK{l})_0$ (resp. $\alpha',
\alpha'' \in \Delta^+(\FRAK{l})$). Let $S$ (resp. $T$) be the set
of indecomposable roots of  $\Delta^+(\FRAK{l}_0)$, (resp.
$\Delta^+(\FRAK{l})).$

\indent
  Let ${\cal O}$ be the
Richardson orbit induced from a Levi factor of $\FRAK{p}_0$. We
say $\FRAK{p}$ is a {\it good parabolic}
if $\dim(\FRAK{g} / \FRAK{p} )_1 = \ell(\cal O)$, (see Section \ref{EM} for notation).
In Section 5 we show that modules induced from a one dimensional module for a good parabolic
 have the least possible multiplicity allowed by the Clifford algebra theory.\\

\subsection{\label{4.2}}
We assume axioms (i)-(iii). For $\lambda \in \FRAK{h}_0^*$ we
define the simple highest weight
 $\FRAK{l}_0$-module $\widehat{L}_S(\lambda)$ as the unique simple
quotient of the Verma module with highest weight $\lambda$ induced
from the Borel subalgebra $ \FRAK{c}(0)$ of $\FRAK{l}_0,$ c.f.
[Ja2,
 5.11].
 The module $\widehat{L}_S(\lambda)$ is finite
dimensional if and only if $\lambda \in P^{++}_{S}$ where
\[ P^{++}_{S} = \{ \lambda \in \FRAK{h}_0^*|(\lambda, \alpha^{\vee})
\in \mathbb{N} \; \mbox{for all} \; \alpha \in S \} . \]

 For $\lambda \in
\FRAK{h}_0^*$ there is a unique graded simple $\FRAK{c}$-module
$V_{\lambda}$ such that $\FRAK{g}^{\alpha}V_{\lambda}= 0$ for all
$\alpha \in \Delta^+(\FRAK{l})$)and $(h - \lambda(h )V_{\lambda} =
0$ for all $h \in \FRAK{h}_0^*$. We remark that if $\FRAK{c}$
involves no classical simple Lie superalgebra of type $Q$, then
$\FRAK{h} = \FRAK{h}_0$ and $dim \; V_{\lambda} = 1$ for all
$\lambda \in \FRAK{h}_0^*.$
 The induced
module $Ind^{\FRAK{l}}_{\FRAK{c}}\;V_{\lambda} $ has a unique
simple graded quotient which we denote by
$\widehat{L}_T(\lambda).$  The conditions for
$\widehat{L}_T(\lambda)$ to be finite dimensional are rather
involved c.f. [K]. However it is easy to show that
 $dim \widehat{L}_T(\lambda) = 1$ if and only if  $\lambda \in
\FRAK{l}^\bot$ where
\[ \FRAK{l}^\bot = \{ \lambda \in \FRAK{h}_0^* | \lambda([\FRAK{l},
\FRAK{l}] \cap \FRAK{h}_0) = 0 \} .\]
 We can regard $\widehat{L}_S(\lambda)$, (resp.
 $\widehat{L}_T(\lambda)$) as a $U({\FRAK{p}}_0)$-module (resp.
 $U({\FRAK{p}})$-module) by allowing  $\FRAK{m}_0^+$,
 (resp. $\FRAK{m}^+$) to act
 trivially and form the induced modules
 \[ M_S (\lambda) = Ind^{\FRAK{g}_0}_{\FRAK{p}_0}\;
 \widehat{L}_S(\lambda), \quad M_T (\lambda) = Ind^{\FRAK{g}}_{\FRAK{p}}\;
 \widehat{L}_T(\lambda) .\]
 By Lemma \ref{IM} \[ d(M_T(\lambda)) = \dim(\FRAK{g}/\FRAK{p})_0
  + d(\widehat{L}_T(\lambda))\]
   and \[e(M_T(\lambda)) =
 2^{(\FRAK{g}/\FRAK{p})_1}e(\widehat{L}_T(\lambda)).\]

\indent To explain the choice of notation:
 $M_S(\lambda)$ conforms to the usage
in [Ja2] while $\widehat{L}_S(\lambda)$ is denoted
$\widehat{L}^S(\lambda)$ in [Ja2].  For $\widehat{L}_T(\lambda)$
and $M_T(\lambda)$ we want something similar which emphasizes the
dependence on $T$ rather than $S$.\\
\\
 \subsection{\label{4.3}}
We can obtain a Lie superalgebra satisfying axioms (i) - (iii) as
follows.
 Suppose that $\FRAK{g}_0$ is reductive with CSA $\FRAK{h}_0$, $V$ is
a $\mathbb{Z}_2$-graded $\FRAK{g}$-module and $V =
\oplus^t_{k=1}V(k)$ where $V(k)$ is a $\mathbb{Z}_2$-graded,
$\FRAK{h}_0$-stable subspace. Set $V(k) = 0$ unless $1 \leq k \leq
t$ and
\[ \FRAK{g}(i) = \{ x \in \FRAK{g}|xV(k) \subseteq V(k-i)\;
\mbox{for all}\; k\}. \]
 In all the examples we consider $\FRAK{g} = \oplus_{i \in \mathbb{Z}} \FRAK{g}(i)$
  satisfies axioms (i) - (iii).\\

  For the rest of this subsection suppose that $\FRAK{g} =
 g\ell(m,n)$ and that ${\bf I} = {\bf I_1} \cup {\bf I_2}$ is as in section \ref{3.1}.
 Set $V = span\{e_i | i \in {\bf I} \}$, the
 natural $\FRAK{g}$-module.
 Consider the function $\sigma: {\bf I} \longrightarrow \{1, \ldots,
 t \}$ defined by $e_i \in V(\sigma(i))$ for all $i$.  Then $e_{ij}
 \in \FRAK{g}(\sigma (j) - \sigma (i))$.
 \indent
 We assume the Borel subalgebra
 $\mathbf{b}$ of upper triangular matrices in $\FRAK{g}_0$ is a
 subalgebra of $\FRAK{p}$.  In terms of $\sigma$ this means that
 if $i < j$ and either $j \leq m$ or $n + 1 \leq i$ we have
 $\sigma (i) \leq \sigma(j)$.
 \indent
 For $1 \leq k \leq t$ set
 \[ \Lambda_k = \{ i \in {\bf I}|\sigma(i) = k \} \]
 and
 \[ r_k = | \Lambda_k \cap {\bf I_1}|, \quad s_k = | \Lambda_k \cap {\bf I_2}| . \]
 We can rearrange the sequences ${\bf{r}} = (r_1, \ldots, r_t)$ and ${\bf{s}} = (s_1,
 \ldots, s_t)$ to obtain partitions $\mu' \in \mathbf{P}(m),
 \nu' \in {\mathbf{P}}(n)$.  Note that the sequences $\bf{r}, \bf{s}$ determine the
 subspaces $V(k).$  Also $ \FRAK{l} = \FRAK{g}(0) \cong
 \oplus^t_{i=1}g\ell(r_i,s_i)$.  It follows that
\[ \dim(\FRAK{g}/ \FRAK{p})_1 = mn - \sum^t_{i=1} r_i s_i . \]
  \begin{lemma} \label{RO}
 (a) The Richardson orbit induced from $\FRAK{l}_0$ is
     ${\cal O}_{\mu,\nu}$.\\
 (b) $\FRAK{p}$ is a good parabolic if and only if there is a permutation
     $\eta$ of $\{ 1, \ldots, t \}$ such that $\mu'_i = r_{\eta(i)}$ and
     $\nu'_i = s_{\eta(i)}$ for $1 \leq i \leq t.$
 \end{lemma}
\noindent{\bf Proof.} (a) follows from [CM, Theorem 7.2.3].\\
 (b) By Theorem \ref{3.3} $\FRAK{p}$ is good if and only if $\sum \mu'_i
 \nu'_i = \sum r_i s_i$.  \\
 We can assume that $\mu'_i = r_i$ for all $i$.  Suppose that $r_j
 > r_{j+1}$ but $s_j < s_{j+1}$ and define $s'_j = s_{j+1}$
 $s'_{j+1} = s_j$ and $s'_i = s_i$ for $i \neq j, j+1$.  Then
 $\sum r_i s'_i > \sum r_i s_i$.  The result follows from this
 observation.\\

We define $\epsilon_i \in \FRAK{h}_0^*$
 so that $\epsilon_i(x)$ is the $i^{th}$ diagonal entry of $x$.
 We take $\Delta^+ = \{ \epsilon_i - \epsilon_j | i < j \}$.  For
 this choice of $\Delta$ we have $\dim \widehat{L}_T(\lambda) <
 \infty \;$ if and only if $\dim\widehat{L}_S (\lambda) < \infty.$  Note that
 $\FRAK{p}$ need not contain the distinguished Borel subalgebra of
 $\FRAK{g}$ as the following examples show.\\
 \\
 {\bf Example.}  Let $(m,n) = (4,3)$ and define $\sigma$ by
 \begin{eqnarray*}
\sigma(1) = \sigma(2) = \sigma(5) & = & 1, \\
 \sigma(3) = \sigma(6) = \sigma(7) & = &2,\\
                         \sigma(4) & = & 3.
 \end{eqnarray*}
 Then \[ \Lambda_1 = \{1, 2, 5\}, \Lambda_2 = \{3, 6, 7\}, \Lambda_3 = \{{4}\}\]
 so ${\bf{r}} = (2, 1, 1), {\bf{s}} = (1, 2, 0).$
Also $S = \{ \epsilon_1 - \epsilon_2, \epsilon_6 - \epsilon_7 \}$
and
\[ \FRAK{l} \cong g \ell (2,1) \oplus g \ell (1,2) \oplus g \ell (1,0)
.\]
 In this case $\FRAK{p}$ is not a good parabolic.

\indent
If we arrange instead that

 \[ \Lambda_1 = \{1, 2, 5,6\}, \Lambda_2 = \{3, 7\}, \Lambda_3 = \{{4}\}\]
 then

\[ \FRAK{l} \cong g \ell (2,1) \oplus g \ell (2,1) \oplus g \ell (1,0)
.\]
 In this case $\FRAK{p}$ is a good parabolic.\\

\subsection {\label{4.4}}
Now let $V, {\FRAK{g}}$ and $\FRAK{g}'$ be as in section
\ref{3.11} and set $\overline{{\FRAK{g}}} = \FRAK{g}'/\FRAK{z}$,
the simple Lie superalgebra of type $Q(n-1)$ . We show how to
associate a good parabolic in $\FRAK{g}$ and $\overline{\FRAK{g}}$
to most nilpotent orbits.
 Suppose that $V =\oplus^t_{k=1}V(k)$ where $V(k)$ is a $\mathbb{Z}_2$-graded
 subspace of $V$
stable under $\FRAK{h}_0$ and $\psi.$ Set $r_k = dim V(k)_0$ and
rearrange the sequence $\mathbf{r}  = (r_1, \ldots , r_t)$ to
obtain a partition $\mu' \in {\bf P} (n).$  The grading on $
\FRAK{g}$ defined in section \ref{4.3} induces a grading on
$\FRAK{g}'$ and
$\overline{\FRAK{g}}.$\\

Let
 $\FRAK{l}_{\mathbf{r}}$ be the block diagonal subalgebra of $ g\ell(n)$
 with diagonal entries of size $r_1, \ldots , r_t$ and set
 $\FRAK{l}'_{\bf r} = \{ x \in \FRAK{l}_{\bf r} | trace(x)= 0 \}.$
 By [CM, Theorem 7.2.3] the Richardson orbit in
$g \ell(n)$ (resp.  $s \ell(n)$ ) induced from $\FRAK{l}_{\bf r}$
 (resp $\FRAK{l}'_{\bf r})$ is ${\cal O}_\mu.$ Also
$\FRAK{g}(0)$ consists  of all matrices of the form
\[ \left[ \begin{array}{cc}
    a & b \\
    b & a
    \end{array}\right] \]
    with $a,b \in \FRAK{l}_{\bf r}$, while $\overline{\FRAK{g}}(0)$
    consists of the images mod $\FRAK{z}$ of matrices of this form
with    $a,b \in \FRAK{l}'_{\bf r}.$ Set $\FRAK{p} = \oplus_{i
\geq 0}   \FRAK{g}(i)$ and $\overline{\FRAK{p}} = \oplus_{i \geq
0} \overline{\FRAK{g}}(i).$  Then $dim (\FRAK{g}/\FRAK{p})_1 = n^2
- \sum_i(\mu'_i)^2.$ Thus from Theorem \ref{3.12} we get the
following result.
 \begin{lemma}\label{Q}
 (a) $\FRAK{p}$ is a good parabolic in $\FRAK{g}$.\\
 (b) If some part of $\mu$ is odd then
 $\overline{\FRAK{p}}$ is a good parabolic in $\overline{\FRAK{g}}$.
 \end{lemma}
\section{Induced Modules and Primitive Ideals}
\subsection {\label{5.1}}
The connection between the Clifford algebras $C \ell_q$ and
modules of low multiplicity is based on the following result.

\begin{lemma}\label{tfree}
Let $N$ be a nonzero finitely generated graded $gr
U(\FRAK{g})$-module such that $q = ann_{S(\FRAK{g}_0)}N$ is prime
and $N$ is torsion free as a $S(\FRAK{g}_0)/q$-module.  If
$\mathcal{V}$ is the closed set in $\FRAK{g}^*_0$ defined by $q$
then $d(N) = d(S(\FRAK{g}_0)/q)$ and $e(N) \geq
2^{\ell(\mathcal{V})}e(\mathcal{V}).$  Furthermore if $e(N) =
2^{\ell(\mathcal{V})}e(\mathcal{V})$ then $C_q$ is split.
\end{lemma}

  {\bf Proof.} Clearly $d(N) \leq d(S(\FRAK{g}_0)/q)$.  Let ${\cal
C} = {\cal C}(q)$ so that $N_{\cal C}$ is a $(gr
U(\FRAK{g})/q)_{\cal C}$-module. There is a factor module of
$N_{\cal C}$ which is a simple module over $C_q = (gr
U(\FRAK{g})/\pi(q))_{\cal C}$.  By [GW, Theorem 9.17 (a)] this
factor has the form $\overline{N}_{\cal C}$ for some $gr
U(\FRAK{g})/\pi(q)$ factor module $\overline{N}$ of $N$. Hence by
the remarks in Section \ref{DM} and Lemma \ref{CA}
\[ d(N) \geq d(\overline{N}) = d(S(\FRAK{g}_0)/q) \] and
\[ e(N) \geq e(\overline{N}) =  2^{\ell(\mathcal{V})}e(\mathcal{V})  . \]
The last statement follows from Lemma \ref{CA}. \vspace*{0.5cm}

\subsection {\label{5.1.1}}
To apply Lemma \ref{5.1} let $M$ be a finitely generated
$U(\FRAK{g})$-module. We equip $M$ with a good filtration and
consider an affiliated series
\[ 0 = N_0 \subset N_1 \subset \ldots \subset N_k = N \]
for the graded module $N = gr M$.  Let $p_1, \ldots, p_k$ be the
affiliated primes of this series and $q_i = \pi^{-1}(p_i)$.  By
[GW, Prop. 2.13] each factor $N_i/N_{i+1}$ is torsion-free as a
$gr U(\FRAK{g})/p_i$-module and hence also as a
$S(\FRAK{g}_0)/q_i$-module.  Thus
\[ e(M) = e(N) = \sum e(N_i/N_{i+1}) \geq \sum 2^{\ell(q_i)}e(S(\FRAK{g}_0)/q_i) \]
where both sums are taken over all indices $i$ such that
$d(N_i/N_{i+1}) = d(N)$.
\\

 By [GW, Proposition 2.14], any prime ideal of which is minimal over
 ann N is equal to one of the $p_i$ and it follows easily that
$ \sqrt{ann_{S(\FRAK{g}_0)} N} = q_1 \cap q_2 \ldots \cap q_k.$
 The closed subset of
  $Spec S(\FRAK{g}_0)$ defined by $ann_{S(\FRAK{g}_0)} N$ is called the {\it associated variety} of $M$.
This definition is independent of the choice of good filtration
[Ja2, 17.2]. These considerations motivate the study of modules
$M$ whose associated variety has a unique component $\mathcal{V}$
with dimension equal to  $d(M)$ and such that $e(M) =
2^{\ell(\mathcal{V})}e(\mathcal{V}).$
\\
For primitive factors $U(\FRAK{g})/P$ the Goldie rank
$rank(U(\FRAK{g})/P)$ is a more important invariant than
$e(U(\FRAK{g})/P)$ so we should try to find primitives $P$ such
that $rank (U(\FRAK{g})/P) \leq 2^{\ell(q)}.$ where $q = \sqrt{ gr
P} \cap S(\FRAK{g}_0).$

\subsection {\label{5.2}}
 For the remainder of the paper we assume that conditions (i) -
(iii) of Section \ref{4.1} hold.
  \begin{theorem}
  Suppose that
  $\dim \widehat{L}_{T}(\lambda) < \infty$
and that $\dim(\FRAK{g}/\FRAK{p})_1 = c$.
  Then $M_T(\lambda)$ has a
filtration by $\FRAK{g}_0$-submodules
\[0 = M_0 \subset M_1 \subset \ldots \subset M_k = M_T(\lambda)\]
such that for $i = 1,\ldots, k$
\[M_i/M_{i-1} \cong M_S(\lambda_i)\]
 for certain $\lambda_i \in P^{++}_{S}$ and
 \begin{equation}\label{DIM}
 \sum^k_{i=1}  \dim \widehat{L}_{S}(\lambda_i) = 2^c \dim
\widehat{L}_T(\lambda).
  \end{equation}
 \end{theorem}
{\bf Proof.}  To simplify notation set $M = M_T(\lambda)$.
 We extend the grading on $\FRAK{g}$ defined in section 4.1 to $U(\FRAK{g})$ and
  $\wedge \FRAK{m}_1$.  Note that $\FRAK{m}_1$ is an
  $\FRAK{l}_0$-module.  Antisymmetrization gives an injective map of
  $\FRAK{l}_0$-modules
   \[ \wedge \FRAK{m}_1 \longrightarrow U(\FRAK{m}) \]
   and we identify $\wedge \FRAK{m}_1$ with its image.  Then
    \[ U(\FRAK{m}) = U(\FRAK{m}_0) \otimes \wedge \FRAK{m}_1 . \]
    It is easy to see that the extended grading satisfies
    \begin{equation}\label{easy}
    [\FRAK{m}^+_0(j), (\wedge \FRAK{m}_1)(-i)] \quad \subseteq \quad \oplus_{r,s,t}
    U(\FRAK{m}_0)(-r)\otimes(\wedge \FRAK{m}_1)(-s) \otimes
    U(\FRAK{p})(t)
    \end{equation}
    for all $i, j > 0$, where the sum is over all $r,s,t \geq 0$
    such that $j-i = t - r - s$. Furthermore since $[\FRAK{m}^+_0(j), \FRAK{m}_1(-k)]
    \subseteq \FRAK{g}(j-k)$ we can restrict the sum on the right to terms with
     $t < j$.  In particular, each
    summand satisfies $s < i$.\\

    \indent Now for $i \geq 0,$ set $L'_i = (\wedge
    \FRAK{m}_1) (-i) \otimes \widehat{L}_T(\lambda) \subseteq M,$
    and
    define inductively $M_0 = 0,$
    and $M'_i = U(\FRAK{g}_0)L'_i + M'_{i-1}$ .  This process
    terminates when $L'_N = (\wedge
    \FRAK{m}_1)\otimes \widehat{L}_T(\lambda)$
    and $M'_N = M$ for some $N.$  Note that each $L'_i$ is an
    $\FRAK{l}_0$-module.
    \indent
    Since $\widehat{L}_T(\lambda)$ is a $U(\FRAK{p})$-module with
    $\FRAK{m}^+_0 \widehat{L}_T(\lambda) = 0$ it follows from
    equation (\ref{easy}) that $\FRAK{m}^+_0 L'_i \subseteq
    M'_{i-1}$.

We refine the series $ 0 = L'_0 \subset L'_1 \subset \ldots
\subset L'_N = \wedge \FRAK{m}_1 \otimes
 \widehat{L}_T(\lambda)$ to a composition series
 \[ 0 = L_0 \subset L_1 \subset \ldots \subset L_k = \wedge \FRAK{m}_1 \otimes
 \widehat{L}_T(\lambda)  \]
of $\wedge \FRAK{m}_1 \otimes
 \widehat{L}_T(\lambda)$ as an $\FRAK{l}_0$-module and define
  $M_i = U(\FRAK{g}_0) L_i + M_{i-1}$.
  Since each $L_i/L_{i-1}$ is  finite
 dimensional it follows that $L_i/L_{i-1} \cong  \widehat{L}_S(\lambda_i)$ for
 $\lambda_i \in P^{++}_{S}.$

    Also for each $i$ we have
    $L'_{j-1} \subseteq L_{i-1} \subset L_i  \subseteq L'_j$
    for some $j$
    and hence
     $\FRAK{m}^+_0 L_i \subseteq \FRAK{m}^+_0 L'_j \subseteq M'_{j-1} \subseteq M_{j-1}$.
      Thus
     \[ \overline{L}_i = (L_i + M_{i-1})/M_{i-1} \]
    is a $U(\FRAK{p}_0)$-module and
    $M_i/M_{i-1} = U(\FRAK{g}_0)\overline{L}_i$.  Hence
    $M_i/M_{i-1}$ is a homomorphic image of
    $Ind^{\FRAK{g}_0}_{\FRAK{p}_0} \; \overline{L}_i$ and  $\overline{L}_i$
    is a homomorphic image of $L_i/L_{i-1}.$ It follows that
 \begin{equation} \label{z1}
    [M_i/M_{i-1}] \leq [Ind^{\FRAK{g}_0}_{\FRAK{p}_0} \; L_i/L_{i-1}].
  \end{equation}
 Therefore
    \[  [M ]= \sum^k_{i=1}  [M_i/M_{i-1}]
    \leq \sum^k_{i=1}  [Ind^{\FRAK{g}_0}_{\FRAK{p}_0} \; L_i/L_{i-1}]
    = [M] \]
where the last equality is obtained by comparing characters using
the PBW theorem.  Thus equality holds in (\ref{z1}) and it follows
that $M_i /M_{i-1} \cong M_S(\lambda_i).$
\\
\\
{\bf Remark.}  From the proof we see that as an
$\FRAK{l}_0$-module
 \[ \oplus_i \widehat{L}_S(\lambda_i) \cong \wedge \FRAK{m}_1 \otimes
 \widehat{L}_T (\lambda) . \]
 With this additional information the theorem generalizes [M1,
 Theorem 3.2].\\
 \\
\subsection {\label{5.3}}
The next result is an analog of [Ja2,15.5(a)].
\begin{corollary}
 If $\widehat{L}_T(\lambda)$ is finite dimensional then  $M_{T}(\lambda)$
  is a homogeneous
  $U(\FRAK{g})$-module.
 \end{corollary}
{\bf Proof.}  Let $N$ be a nonzero submodule of $M_{T}(\lambda)$
and choose $i$ minimal such that $N \cap M_i \neq 0$.  Then $N
\cap M_i$ is isomorphic to a nonzero submodule of $M_S(\lambda_i)$
which is a homogeneous $U(\FRAK{g}_0)$-module by [Ja2, Satz
15.5(a)]. Hence
\begin{eqnarray*}
d(M_S(\lambda_i)) &=& d(N \cap M_i) \\
                  & \leq & d(N) \leq d(M_{T}(\lambda)) .
\end{eqnarray*}
The result follows since $d(M_S(\lambda_i)) = d(M_T(\lambda)) =
\dim(\FRAK{g}/\FRAK{p})_0$.\\

\subsection {\label{5.4}} The following result is an analog of
[Ja2, 17.16].
\begin{lemma}
The associated variety $V(gr ann_{U(\FRAK{g}_0)} M_T(\lambda))$ is
the closure of the Richardson orbit induced from a Levi factor of
$\FRAK{p}_0.$
\end{lemma}
{\bf Proof.} Consider the series $M_0 \subset M_1 \subset \ldots
\subset M_k = M_{T}(\lambda)$ of Theorem \ref{5.2} and set
$I_S(\lambda_i) = ann_{U(\FRAK{g}_0)} M_S(\lambda_i).$  Then
\[ I_S(\lambda_1) \ldots I_S(\lambda_k) \subseteq ann_{U(\FRAK{g}_0)}
 M_{T}(\lambda)\subseteq I_S(\lambda_i) \]
so that
\[ \prod gr I_S(\lambda_i)  \subseteq gr ann_{U(\FRAK{g}_0)}
    M_{T}(\lambda)\subseteq gr I_S(\lambda_i) \] for all $i$.
  On the other hand by [Ja2, 17.16] $V(gr ann_{U(\FRAK{g}_0)} M_S(\lambda_i)) =
   G \FRAK{m}_0$ for
  all $i$ so the result follows.\\
\\
\subsection {\label{5.5}}  For the proof of Theorem \ref{5.6} we
need a good filtration on  $M_T(\lambda)$ with special properties.
\begin{lemma}
   If $M = M_T(\lambda)$ and $q = S(\FRAK{g}_0)\FRAK{p}_0$,
there is a good filtration on $M$ such that
   $ann_{S(\FRAK{g}_0)}(gr M) = q$ and $gr M$ is a torsion free
   $S(\FRAK{g}_0)/q$-module.
 \end{lemma}
{\bf Proof.}   Let $U_n(\FRAK{m}) = U_n(\FRAK{g}) \cap
U(\FRAK{m})$
 and
\[ M_n = U_n(\FRAK{m}) \otimes \widehat{L}_T(\lambda) \]
Since $\FRAK{p} \widehat{L}_T(\lambda) \subseteq
\widehat{L}_T(\lambda)$, an easy induction shows that $\FRAK{g}_0
M_n \subseteq M_{n+2}$ and $\FRAK{g}_1M_n \subseteq M_{n+1}$, so
$\{M_n\}$ is a filtration of $M$ as a $U(\FRAK{g})$-module.

Similarly we have $\FRAK{p}_0 M_n \subseteq M_n$ which implies
$S(\FRAK{g}_0)\FRAK{p}_0 \subseteq ann_{S(\FRAK{g}_0)} gr M$.  On
the other hand $gr M \simeq gr U(\FRAK{m}) \otimes
\widehat{L}_T(\lambda)$ is a free $S(\FRAK{m}_0)$-module. The
result follows from this.\\

\subsection {\label{5.6}}
Part (a) of the next result is an analog of [Ja2, Satz 15.5b)].

\begin{theorem}
 Suppose that $\FRAK{p}$ is a good parabolic in $\FRAK{g}$ and that
$\dim \widehat{L}_T(\lambda) = 1$.  Set $M = M_T(\lambda)$ and $q
= S(\FRAK{g}_0)\FRAK{p}_0$ . Then\\
 (a) $M$ is a critical $U(\FRAK{g})$-module with $e(M) = 2^{\ell(q)}$\\
 (b) $ann_{U(\FRAK{g})} M$ is a primitive ideal.
 \end{theorem}
 {\bf Proof.}  (a) Consider a good filtration on $M$ as in Lemma
 \ref{5.5}.  If $M'$ is a nonzero submodule of $M$ then $N' = gr M'$ is
 a nonzero
 submodule of $N = gr M$, and we have $d(N') = d(M')$ and $e(N') =
 e(M')$.
\indent Let $q'$ be the prime ideal of $S(\FRAK{g}_0)$ defining
the Richardson orbit induced from $\FRAK{l}_0$ and $q =
S(\FRAK{g}_0)\FRAK{p}_0$.  Since $\FRAK{p}$ is a good parabolic
\[ \dim(\FRAK{g}/\FRAK{p})_1 = \ell(q). \]

\indent
 By Lemma \ref{5.5} $N'$ is torsionfree, so by Lemma \ref{5.1} we have
 \[ \dim (\FRAK{g}/\FRAK{p})_0 \leq d(N') \leq d(M) = \dim(\FRAK{g}/\FRAK{p})_0 \]
 and
 \[ 2^{\dim(\FRAK{g}/\FRAK{p})_1} \leq e(N') \leq e(M) =
 2^{\dim(\FRAK{g}/\FRAK{p})_1} . \]
 Thus equality holds in both cases and this proves the result.\\
 \\
 (b)  Note that $M$ has finite length, so the arguments in [Ja2,
      8.14-8.15] show that $soc M$ is simple and $ann_{U(\FRAK{g})} M
      = ann_{U(\FRAK{g})} soc M$.\\
\\
\subsection {\label{5.7}}
 If $M,N$ are $U(\FRAK{g}_0)$-modules
we set as in [Ja2]
 \[ {\cal L}(M,N) = \{ \phi \in Hom_{\mathbb{C}}(M,N)|\dim
 U(\FRAK{g}_0)\phi < \infty\} .\]
 Then ${\cal L}(M,N)$ is a $U(\FRAK{g}_0)$-bimodule, and if $M,N$
 are actually $U(\FRAK{g})$-modules then ${\cal L}(M,N)$ is a
 $U(\FRAK{g})$-bimodule.

 For $X$ a  $U(\FRAK{g}_0)$-bimodule we write $R ann X$
 for the annihilator of $X$
 as a right  $U(\FRAK{g}_0)$-module. If $\wedge$ is a coset of the integral weight lattice of
 $\FRAK{g}_0$ in $\FRAK{h}^*_0$ the set $\wedge^{++}$ is defined
 as in [Ja2, 2.5].
 \begin{lemma}\label{prim}
  If $\lambda \in \wedge^{++}, \widehat{L}_T(\mu)$ is finite dimensional and
  $M = Ind^{\FRAK{g}}_{\FRAK{p}}\; \widehat{L}_T(\mu)$ then
  $R ann_{U(\FRAK{g}_0)}{\cal L}(M(\lambda), M)$ is a primitive
  ideal of $U(\FRAK{g}_0)$.
   \end{lemma}
   {\bf Proof.}
   There is a surjective map of
   $U(\FRAK{g}_0)$-modules
   \[ M' = U(\FRAK{g}) \otimes_{U(\FRAK{p}_0)} \widehat{L}_T(\mu)
   \longrightarrow M. \]
 Since the finite dimensional module $\widehat{L}_T(\mu)$ is
 semisimple as a $\FRAK{l}_0$-module, we
    can write \[ \widehat{L}_{T}(\mu) \cong \oplus
   \widehat{L}_{S}(\mu_i) . \]
   with $\mu_i \in \wedge \cap P_S^{++}.$
   Thus as a $\FRAK{g}_0$-module
   \[ M' \cong \bigoplus_i U(\FRAK{g})\otimes_{U(\FRAK{g}_0)} U(\FRAK{g}_0)
\otimes_{U(\FRAK{p}_0)} \widehat{L}_{S}(\mu_i) \]
\[ =    \bigoplus_i U(\FRAK{g}) \otimes_{U(\FRAK{g}_0)}
   M_S(\mu_i) \]
\[ \cong \bigoplus_i E \otimes_{U(\FRAK{g}_0)} M_S(\mu_i) \]
 for some finite dimensional $U(\FRAK{g}_0)$-module E.  Since
 $\lambda \in \wedge^{++}$ the functor ${\cal L}(M(\lambda), \underline{\;\;\;}  )$
 is exact on the category ${\cal O}$, see [Ja2, Lemma 4.8 and 6.9 (9)].  Hence
 $X' = {\cal L}(M(\lambda),M')$ maps onto $X = {\cal
 L}(M(\lambda),M)$ and $R ann X \subseteq R ann X'$.  Similarly
 since $M'' = M_S (\mu) \subseteq M_T(\mu)$ we have $X'' = {\cal
 L}(M(\lambda), M'') \subseteq {\cal L}(M(\lambda),M)$ and so $R
 ann X'' \supseteq R ann X$.   Finally [Ja2, $6.8\;(2')$ and Lemma 15.7]
 imply that $R ann X'' = R ann X'$ is a primitive ideal in
 $U(\FRAK{g}_0)$.

\subsection {\label{5.8}}  By Corollary \ref{5.3} and Lemma \ref{5.7} the hypotheses
 of [Ja2, Satz 12.3] are satisfied.  We apply this below.
 \begin{theorem}\label{GR}
 If $\dim \widehat{L}_T (\lambda) < \infty$, $M = M_{T}(\lambda)$ and $c =
 \dim(\FRAK{g} / \FRAK{p})_1$ then ${\cal L}(M,M)$ is prime Noetherian with
 Goldie rank $2^c \dim \widehat{L}_T (\lambda)$.
 \end{theorem}
 {\bf Proof.}   By [Ja2, Satz 12.3 (a), (c) ]
  ${\cal L}(M,M)$
  is prime Noetherian and
\[ \mbox{rank}\;  {\cal L}(M,M) = \sum_L [M:L]\; \mbox{rank}\; {\cal L}(L,L) \]
 where the sum runs over composition factors $L$ of $M$ as a $U(\FRAK{g}_0)$-module such that
 $d(L) = d(M)$, and [M : L] is the multiplicity of $L$ in $M$.
 Now if
 \[ 0 = M_0 \subset M_1 \subset \ldots \subset M_k = M \]
 is the series given in Theorem \ref{5.2} then
  \[ [M : L] = \sum^k_{i=1} [M_S (\lambda_i): L ] . \]
Using [Ja2, 15.8], then [Ja2, 15.21 (2)] and finally equation
(\ref{DIM}) we obtain
\begin{eqnarray*}
\mbox{rank}\; {\cal L}(M,M) & = & \sum^k_{i=1} \mbox{rank}\; {\cal
L}(M_S(\lambda_i), M_S(\lambda_i)) \\
  & = & \sum^k_{i=1} \dim \widehat{L}_S(\lambda_i) = 2^c dim \widehat{L}_T(\lambda).
 \end{eqnarray*}
\subsection {\label{5.10}} In the final result of this subsection we assume that $\FRAK{g} = \oplus_{i
\in \mathbf{z}}\FRAK{g}(i)$ is a graded Lie superalgebra as in
section 4.1 and set $\FRAK{p} = \oplus_{i \geq 0}\FRAK{g}(i),
\FRAK{l} = \FRAK{g}(0)$.  Suppose that $\lambda \in \FRAK{l}^\bot$
and set $M = M_T(\lambda), q =
S(\FRAK{g}_0)\FRAK{p}_0$.\\
\\
{\bf Corollary.}  {\it Suppose that $\FRAK{p}$ is a good parabolic
in $\FRAK{g}$.  Then
 $e(M) = 2^{\ell(q)},$ $ann_{U(\FRAK{g})}M$
  is a prime ideal of $U(\FRAK{g})$ and
${\cal L}(M,M)$ is a primitive ring with Goldie rank
$2^{\ell(q)}$.  Furthermore the Clifford algebra
$C_q$ is split. }\\
{\bf Proof.} Since dim $\widehat{L}_T(\lambda) = 1$, the
statements about ${\cal L}(M,M)$ follow from Theorem \ref{5.8},
while the claims about $M$ and $ann_{U(\FRAK{g})}M$ hold by
Theorem \ref{5.6}. Since
$e(M) = 2^{\ell(q)}$, $C_q$ is split by Lemma  \ref{CA} .\\

Observe that $U(\FRAK{g})/ann_{U(\FRAK{g})}M$ embeds in
$\mathcal{L}(M,M)$.    It follows from a version of the additivity
principle [GW, Corollary 7.26] that $\overline{U} =
U(\FRAK{g})/ann_{U(\FRAK{g})}M$ has Goldie rank at most
$2^{\ell(q)}$.  However the Goldie rank of $\overline{U}$ can be
strictly less than $2^{\ell(q)}$. For example suppose that
$\FRAK{g} = osp(1,2)$ and let $M$ be a Verma module.  The
associated variety of $M$ in $\FRAK{g}_0$ is the nilpotent cone
$\mathcal{N}$ and we have $k(\mathcal{N}) = \ell(\mathcal{N}) =
1,$ by Theorem \ref{3.8}.  By the Corollary the Goldie rank of
$\mathcal{L}(M,M)$ is 2, but for an appropriate choice of $M,
\overline{U}$ is isomorphic to the first Weyl
algebra, which has Goldie rank 1,  see [P].\\

\section{Nilpotent Orbits in Lie superalgebras}
\subsection {\label{km}}
Suppose that $\FRAK {g}$ is a Lie superalgebra such that ${\FRAK
g}_0$ is reductive, and that there is a non-degenerate even
invariant bilinear  form $B$ on ${\FRAK g}.$  We use $B$ to
identify ${\mathfrak g}_0^*$ with ${\mathfrak g}_0.$  For $x \in
{\mathfrak g}_0,$  let ${\mathfrak g}^x$ be the centralizer of $x$
in $ {\mathfrak g}$.
 \begin{lemma} We have
\[ k(m_x) =
dim {\mathfrak g}_1 - dim {\mathfrak g}^x_1. \]
   \end{lemma}
   {\bf Proof.}
If $u \in {\mathfrak g}_1,$ then $u \in {\mathfrak g}^x$ if and
only if
     $0 = B([x, u],w)=B(x,[u,w])$ for all $w \in \mathfrak g_1.$
This holds if and only if $u$ is in the radical of the
$\mathbb{C}$-valued bilinear form on ${\mathfrak g}_1$ whose
matrix is obtained by reducing $M({\mathfrak g})$ mod $m_x.$ The
result follows since $k(m_x)$ is the rank of this bilinear form.\\

 \noindent {\bf Remarks.} If $ad \; x$ is nilpotent, then
the value of $k(m_x)$ is given by the formulas in Section 3. When
$\FRAK{g} = g \ell(m,n),$ and $ x \in {\mathfrak g}_0,$  we can
compute $dim {\mathfrak g}^x_1 $ directly as follows.  If \[ x =
\left[
\begin{array}{cc} J_\mu & 0\\0 & J_\nu \end{array}\right] ,y = \left[
\begin{array}{cc} 0 & C\\D & 0
\end{array}\right], \]
we have $y \in {\mathfrak g}^x_1 $ if and only if $J_\mu C = C
J_\nu$ and $D J_\mu =  J_\nu D.$ Let $U =  \{ C \in  M_{m,n}|
J_\mu C = C J_\nu \}.$ Then $U$ is the space of highest weight
vectors in $Hom( \oplus_{i \geq 1} L(\mu_i), \oplus_{i \geq 1}
L(\nu_i)) .$ Hence as in the proof of Lemma \ref{Part}, we have
$dim \; U = \sum_i \mu'_i \nu'_i .$ This easily gives $ dim
{\mathfrak g}^x_1 = 2 \sum_i \mu'_i \nu'_i $.

\subsection {\label{6.2}}
Consider the action of an algebraic group $K$ on its Lie algebra
${\mathfrak k}$ by the adjoint representation $Ad:K
\longrightarrow GL({\mathfrak k}).$  It follows from [H,Theorem
10.4] that the tangent space to the orbit at $x \in \mathfrak k$
is given by $ T_x(K \cdot x) = {\mathfrak k}/{\mathfrak k}^x.$ We
prove a parallel result for certain Lie superalgabras. Before we
can state it however, we need to review some notions concerning
the functor of points and Lie supergroups,
see [DG],[Ja3] and [Ma].\\

The category of supercommutative $\mathbb{C}$-algebras will be
 denoted {\bf Alg} and the category of sets by {\bf Set}.
Whenever we construct a functor $X$ from {\bf Alg} to {\bf Set},
we do so by specifying the value of $X$ on an object  $R$ of {\bf
Alg} in a way which is functorial in $R.$  Hence there is no need
to say anything about the effect of $X$ on morphisms. We call
$X(R)$ the set of $R$-points of $X.$ We say that $X$ is a {\it
subfunctor} of $Y$ if $X(R) \subseteq Y(R)$ for all
supercommutative $R.$ An {\it affine superscheme} $X$ is a
representable functor from {\bf Alg} to {\bf Set}.  Thus there is
a supercommutative $\mathbb{C}$-algebra $\mathcal{O}(X)$ such that
$X(R)  =  Mor_{\bf Alg }(\mathcal{O}(X),R)$ for any
supercommutative algebra $R.$

\subsection{}
 Suppose that $V = V_0 \oplus V_1$ is a $\mathbb{Z}_2$-graded vector
 space and let $V^*$ be the dual vector space.  To specify $V$ as a
 representable functor we need to define $\mathcal{O}(V).$  This
 is done by setting
 \[ \mathcal{O}(V) = S(V^*_0) \otimes \wedge(V^*_1), \]
the tensor product of the symmetric algebra on the vector space
 $V^*_0$
 and the exterior algebra
 on the vector space $V^*_1.$
It is easy to see that for any supercommutative algebra $R$
\[ V(R) = V_0 \otimes R_0 + V_1 \otimes R_1 .\]
It should be clear from the context whether $V$ is to be thought
of as a  $\mathbb{Z}_2$-graded vector space or as a functor.  We
say that an affine superscheme $X$ is a {\it closed subscheme} of
$V$ if $\mathcal{O}(X)$ is a $\mathbb{Z}_2$-graded factor algebra
of $\mathcal{O}(V).$  If ${\mathfrak g}$ is a Lie superalgebra,
then for any supercommutative algebra $R, \; {\mathfrak g}(R)$
becomes a Lie algebra when we set
\[ [u \otimes r, v \otimes s] = [u, v]rs \]
for all $u \otimes r \in {\mathfrak g}_i \otimes R_i, v \otimes s
\in {\mathfrak g}_j \otimes R_j, \; (i,j = 0, 1).$

\subsection{} Let $H$ be a supercommutative Hopf superalgebra with coproduct $\Delta.$
 For $h \in H,$ write
 \[ \Delta (h) = \sum h_1 \otimes h_2 .\]
 Then for any $R \in Ob \; {\bf Alg}$,  $Hom_{\mathbb{C}}(H,R)$ is an algebra
 under the convolution product
\begin{equation} \label{eq1}
 (\phi \cdot \omega)(h) = \sum \phi(h_1)\omega(h_2)
\end{equation}
 for \[\phi,\omega \in Hom_{\mathbb{C}}(H,R).\]
Note that the identity of $Hom_{\mathbb{C}}(H,R)$ is the composite
of the counit $H \longrightarrow \mathbb{C}$ followed by the
inclusion $\mathbb{C} \longrightarrow R$. Also $Mor_{\bf Alg
}(H,R)$ is a subgroup of the group of units of
$Hom_{\mathbb{C}}(H,R).$ The inverse of $\phi \in Mor_{\bf Alg
}(H,R)$ is $\phi^{-1} = \phi \circ \sigma$ where $\sigma$ is the
antipode of $H$.

\subsection{}
Let $V$ be a $\mathbb{Z}_2$-graded vector space and $K = GL(V).$
By choice of a basis we identify $K$ with $GL(m,n).$  For $R \in
Ob \; {\bf Alg}$ the $R$-points of $K$ are matrices over $R.$  We
use the set ${\bf I} = {\bf I_1} \cup {\bf I_2}$  as in section
\ref{3.1} to index the rows and columns of these matrices, as well
as elements of $\FRAK{k} = g \ell(m,n).$ We think of $K$ as the
group scheme represented by the Hopf superalgebra $H =
\mathcal{O}(K)$ which we describe below.

Treating  $\FRAK{k}$ as a  $\mathbb{Z}_2$-graded vector space, the
construction of the previous subsection yields an algebra
$\mathcal{O}(\FRAK{k}).$ It is often convenient to arrange the
generators $x_{ij}, (i,j \in {\bf I})$ of $\mathcal{O}(\FRAK{k})$
in {\it standard matrix format} $[Ma, page 158]$.  This means that
we arrange them in the form
\[ x = (x_{ij}) = \left[ \begin{array}{cc}
x_1 & x_2 \\
 x_3 & x_4 \end{array}\right] \]
 where $x_1$ is the matrix of indeterminates $(x_{ij})_{i,j \in \bf I_1}$ and
the other submatrices are defined similarly.  All entries in the
matrices $x_1, x_4$ are even while those in $x_2, x_3$ are odd. As
an algebra $\mathcal{O}(\FRAK{k})$ is the tensor product of the
polynomial algebra generated by the even entries of $x$ with the
exterior algebra on the vector space spanned by the odd entries of
$x.$

We can make  $\mathcal{O}(\FRAK{k})$ into a bialgebra by defining
the coproduct $\Delta$ and counit $\epsilon$ on the generators
$x_{ij},\;\;
 (i,j \in {\bf I})$ by
\begin{eqnarray*}
\Delta(x_{ij}) & = & \sum_{\ell \in {\bf I}} x_{i\ell} \otimes
x_{\ell
j} \\
 \epsilon(x_{ij}) & = & \delta_{ij}.
\end{eqnarray*}
This implies that the product defined by equation (\ref{eq1})
 is just matrix multiplication.

Note that $d = (\det x_1)(\det x_4)$ is a polynomial in the
central variables $x_{ij}, x_{kl}$ with $i, j \in {\bf I}_1$, $k,l
\in {\bf I}_2$.  Inverting $d$ we obtain the Hopf superalgebra $H
= \mathcal{O}(K)$.  The coproduct and counit for $H$ are uniquely
determined by their counterparts for
$\mathcal{O}(\mathcal{M}_{m,n})$.  The antipode $\sigma$ is
defined on generators $x_{ij}$ symbolically by
\[ \sigma\left( \left[ \begin{array}{cc}
x_1 & x_2 \\
 x_3 & x_4
 \end{array}\right]\right) = \left[ \begin{array}{cc}
  y_1 & -x^{-1}_1 x_2 y_4 \\
-x^{-1}_4 x_3 y_1 & y_4 \end{array} \right] \]
 where
\[ y_1 = (x_1 - x_2 x^{-1}_4 x_3)^{-1}, \quad y_4 = (x_4 - x_3 x^{-1}_1
x_2)^{-1} . \] Thus for example if $i,j \in \bf I_1$, then
$\sigma(x_{ij})$ is the entry in row $i$ and column $j$ of $y_1.$

We say that $G$ is a {\it closed subgroup} of $GL(V)$ if $G$ is a
subfunctor of
 $Mor_{\bf Alg }(H,\;\;)$ of the form $G = Mor_{\bf Alg }(H/I,\;\;)$
 for some Hopf ideal $I$ of $H.$  If this is the case we set
 $\mathcal{O}(G) = H/I.$ We say that the functor $G$ is a {\it (linear) Lie
 supergroup} if it is isomorphic to a closed subgroup of $GL(V)$ for some
 $V.$  A Lie supergroup $G$ {\it acts} on an affine superscheme $X$ if $\mathcal{O}(X)$
 is an $\mathcal{O}(G)$-comodule
 algebra, [Mo,4.1.2]. This means that there is a natural transformation
of functors $G \times X \longrightarrow X$ satisfying the usual
axioms for group actions, see [DG, page 160].

\subsection{}
We define orbits for actions of Lie supergroups and study their
tangent spaces. Suppose that $G$ is a Lie supergroup which acts on
an affine superscheme $X.$  We write $X_{red}$ for the
$\mathbb{C}$-points of $X,$ that is \[ X_{red} = Mor_{\bf Alg
}(\mathcal{O}(X),\mathbb{C}).\] If $x \in X_{red},$ then using the
inclusion $\mathbb{C} \longrightarrow R$ we can regard $x$ as an
element of $X(R)$ for any supercommutative algebra $R.$  It is
therefore meaningful to define
 a subfunctor $G\cdot x$ of $X$
by setting
\[(G\cdot x)(R) = \{ g\cdot x | g \in G(R) \},\]
compare [DG, page 243]. The orbit map $ g \rightarrow g\cdot x$
gives a natural transformation of functors $\mu:G \rightarrow
{G\cdot x} \subseteq X.$ The orbit closure $\overline{G\cdot x}$,
can be described as follows, cf [S1]. By Yoneda's lemma, $\mu$
induces an algebra map $\mu^*:\mathcal{O}(X) \longrightarrow
\mathcal{O}(G),$ and $\overline{G\cdot x}$ is the closed
subfunctor of $X$  defined by the ideal $Ker \mu^*.$

\subsection {\label{TH1}}
For a supercommutative algebra $R,$ the {\it algebra of dual
numbers
 over} $R$ is the algebra $R[\varepsilon]$ obtained from $R$ by
adjoining an even central indeterminate $\varepsilon$ such that
$\varepsilon^2 = 0.$  Suppose that $X$ is a subfunctor of a
$\mathbb{Z}_2$-graded vector
 space $V.$  We define the tangent space,
$T_x(X)$ to $X$ at $x \in X_{red}.$ As a first attempt we consider
the subfunctor of $V$ given by
\[t_x(X)(R) =
\{y \in Hom_{\mathbb{C}}(\mathcal{O}(V),R)| x + y\varepsilon \in
X(R[\varepsilon]) \}. \] That is  $t_x(X)(R)$ is the fiber over
$x$ under the map $X(R[\varepsilon]) \longrightarrow X(R).$
However it is not clear that $t_x(X)$ is a subspace of $V,$
compare the discussion in [EH, VI.1.3] on the tangent space to a
functor. So we define $T_x(X)$ to be the smallest
$\mathbb{Z}_2$-graded subspace of $V$ such that $t_x(X)(R)
\subseteq T_x(X)(R)$  for all supercommutative $R.$
 For $X \subseteq K = GL(V),\; t_x(X)$ and $T_x(X)$ are
defined in similar ways except that $V$ is replaced by ${\mathfrak
k} = g\ell (V).$  This definition works well in the cases of
interest to us which
are as follows.\\

{\bf Case 1.} Suppose that $X$ is a closed subscheme of $V$ and
let ${\bf{ m}}_x$ be the maximal ideal of $\mathcal{O}(X)$
corresponding to $x.$ Using a Taylor expansion centered at $x,$
see [Le,II.2], we can see that $t_x(X) = T_x(X)$ is naturally
isomorphic to $({\bf m}_x/{\bf m}^2_x)^*.$  Note also that if $y
\in Hom_{\mathbb{C}}(\mathcal{O}(V),R),$
 then $y \in T_x(X)(R)$ if and only if $y$ is an $R$-valued
point derivation at $x$ (compare [H, page 38]), that is
\[ y(fg) = x(f)y(g) + y(f)x(g) \]
for all $f, g \in \mathcal{O}(X).$

{\bf Case 2.} If $G$ is a closed subgroup of $K,$ the tangent
space to the identity $1 \in G$ is
\[T_1(G)(R) = \{y \in Hom_{\mathbb{C}}(\mathcal{O}({\mathfrak
k}),R) | 1 + y\varepsilon \in G(R[\varepsilon]) \}. \] An easy
computation, [A, Chapter 8, (6.19)], shows that
$T_1(G)(R)$ is a Lie subalgebra of ${\mathfrak k}(R)$ for any
supercommutative algebra $R.$  Thus $T_1(G)$ is a Lie superalgebra
which we denote
by $Lie(G).$\\

{\bf Case 3.} For $G, K$ as in Case 2, set  ${\mathfrak g}
=Lie(G)$ and ${\mathfrak k} =Lie(K).$ Then $G(R)$ acts by
conjugation on ${\mathfrak g}(R).$ Since the action is functorial
in $R,$ we can say that $G$ acts on ${\mathfrak g}.$ Consider the
orbit of a $\mathbb{C}$-point of $ {\mathfrak g},$ that is of an
element $x \in \mathfrak{g}_0.$ We have
\[t_x(G\cdot x)(R)= \{ y \in Hom_{\mathbb{C}}(\mathcal{O}({\mathfrak k}),R)| x +
y\varepsilon \in ({G\cdot x})(R[\varepsilon]) \}. \]

To compute this, suppose $g_0 \in G(R)$ and $g_1 \in
Hom_{\mathbb{C}}(\mathcal{O}( {\mathfrak k}),R)$ and set $g = g_0
+ g_1\varepsilon.$  Then
 $g^{-1} =  g^{-1}_0  -  g_0^{-1}
 g_1g^{-1}_0\varepsilon.$  Since $G(R) \subseteq  G(R[\varepsilon])$ we have $g
\in (G(R[\varepsilon])$ if and only if $gg^{-1}_0 = 1 +
g_1g^{-1}_0 \varepsilon \in G(R[\varepsilon])$  and this is
equivalent to $z = g_1g^{-1}_0 \in \mathfrak{g}(R).$

We have $gxg^{-1} = x + y\varepsilon$
 if and only if $g_0xg^{-1}_0 =x$ and
$gxg^{-1} = x +  [z,x]\varepsilon .$  The calculations take place
in the algebra $Hom_{\mathbb{C}}(\mathcal{O}( {\mathfrak
k}),R[\varepsilon]).$
 It follows that $x +
y\varepsilon \in ({G\cdot x})(R[\varepsilon])$ if and only if $y =
[z,x] \in [{\mathfrak g}(R),x] = [{\mathfrak g},x](R).$
 Therefore $t_x (G\cdot x) =
[{\mathfrak g},x],$  which is a subspace of ${\mathfrak g}.$  We
have proved the following result.
\begin{theorem} We have
\[T_x (G\cdot
x) =  [{\mathfrak g},x] \] and the  map $z \longrightarrow [z,x]$
gives a natural isomorphism of functors
\[{{\mathfrak g}/{\mathfrak g}^x} \longrightarrow  T_x (G\cdot
x) .\]
\end{theorem}
{\bf Remark.}  We do not know whether, in the situation of
the Theorem, we have $T_x (G\cdot x) = T_x (\overline{G\cdot x}).$\\

 We define
the {\it superdimension} of a $\mathbb{Z}_2$-graded vector
 space $U = U_0 \oplus U_1$ to be
 \[ \underline{dim} \; U = (dim \; U_0 , dim \;U_1) .\]
In the non-super case, if $X$ is an irreducible variety, we have
$dim \; T_x(X) \geq dim \; X$ with equality on a dense subset of
$X,$ [H,Theorem 5.2].  For an orbit $X = G\cdot x$ as above, it
only makes sense to consider the tangent space at a
$\mathbb{C}$-point  $y$ of $X.$  In this case, clearly $G\cdot x =
G\cdot y$ so $T_y(G\cdot x) = T_y(G\cdot y)$ has the same
dimension as $T_x(G\cdot x).$ Hence it is reasonable to define
$\underline{dim} \; G\cdot x$ to be $\underline{dim} \; T_x
(G\cdot x).$

\subsection {\label{6.11}}
Since there are Lie algebras which are not  the Lie algebra of any
algebraic group, see [McR,14.7.4], the question now arises whether
Theorem \ref{TH1} applies to classical simple Lie superalgebras.
This is the case at least in the following examples.\\

{\bf Example 1.}  If $G =  GL(V),$ then clearly $Lie(G)  = g\ell (V).$ \\

{\bf Example 2.}  Let $Ber \in \mathcal{O}(GL(V))$ be the
Berezinian, or superdeterminant, [Ma, Section 3.3].  This is a
grouplike element of $\mathcal{O}(GL(V)).$ We define $SL(V)$ to be
the group scheme represented by the Hopf superalgebra
$\mathcal{O}(GL(V))/(Ber - 1).$ It is well known that if $G =
SL(V)$, then $Lie(G)  = s\ell (V).$    This is
easy to see using our definition of $T_1(G).$\\

{\bf Example 3.}  Let $K = GL(V),$ and ${\mathfrak k} =  g\ell
(V)$ and suppose that $(\;,\;)$ is a homogeneous bilinear form on
$V.$ The Lie superalgebra ${\mathfrak g}$ preserving this form is
defined, see [Sch page 129],  by setting
\[ {\mathfrak g}_a = \{ x \in {\mathfrak k}_a | (xu,v) +
(-1)^{a \overline{u}}(u,xv) = 0 \; \mbox{for all} \; u, v \in V,
\mbox{with} \; deg \; u = \overline{u}\} .\] We extend $(\;,\;)$
to a bilinear form $(\;,\;)_R$ on $V(R)$ by the rule
\[ (u \otimes r, v \otimes s)_R= (u, v)rs \]
for all $u \otimes r \in V_i \otimes R_i, v \otimes s \in V_j
\otimes R_j, \; (i,j = 0, 1)$.  It is easy to show that if $R_1$
is sufficiently large, then
\[ {\mathfrak g}(R) = \{ g \in \mathfrak k(R)| (gu,v)_R + (u,gv)_R = 0 \; \mbox{for all} \; u, v \in
V(R) \}. \] That is  ${\mathfrak g}(R)$ is the Lie subalgebra of
${\mathfrak k}(R)$ preserving the form $(\;,\;)_R$. On the other
hand the Lie supergroup $G$ preserving $(\;,\;)$ is the functor
defined by
\[ G(R) = \{ g \in K(R)| (gu,gv)_R = (u,v)_R \; \mbox{for all} \; u, v \in
V(R) \} .\] A simple computation shows that $Lie(G) = {\mathfrak
g}.$

\begin{theorem} Let ${\mathfrak g} =  s\ell (m,n), (m
\neq n), g\ell (m,n),$ or $osp(m,n)$ and let $G$ be the Lie
supergroup with $Lie(G) =  {\mathfrak g}$ defined above. Then if
$x \in {\mathfrak g}_0$ we have \[ \underline{dim}  \; G\cdot x =
(dim \; G_0\cdot x, k(m_x)) .\]
\end{theorem}
   {\bf Proof.}
In these cases there is a non-degenerate even invariant bilinear
form on ${\FRAK g}.$ Hence by Lemma \ref{km} we have $k(m_x) = dim
{\mathfrak g}_1 - dim {\mathfrak g}^x_1 .$ But by Theorem
\ref{TH1} we have $dim \; T_x (G\cdot x)_1 =  dim \; {\mathfrak
g}_1/{\mathfrak g}^x_1.$  This proves the statement about $dim \;
T_x(G\cdot x)_1,$ and the claim about $dim \; T_x(G\cdot x)_0$
follows similarly.\\

We remark that the values of $dim \; G_0\cdot x$ for nilpotent
orbits in classical Lie algebras are given in [CM, Corollary
6.1.4].

  \pagebreak
\begin{center}
{\large{\bf{References}}}
\end{center}
\begin{itemize}
\item[{[A]}] M. Artin,  Algebra. Prentice Hall, Inc., Englewood
Cliffs, NJ, 1991. \item[{[B]}] A.D. Bell, A criterion for
primeness of enveloping algebras of Lie superalgebras. J. Pure
Appl. Algebra  69  (1990),   111--120. \item[{[Be]}] E. Benlolo,
Sur la quantification de certaines variet\'{e}s orbitales, Bull.
Sci. Math. 118 (1994), 225-243. \item[{[BB]}] W. Borho and J.-L.
Brylinski, Differential operators on homogeneous spaces I, III,
Invent. Math. 69 (1982), 437-476 and 80 (1985), 1-68.
\item[{[CM]}] D. H. Collingwood and W. M. McGovern, Nilpotent
Orbits in Semisimple Lie Algebras, Van Nostrand Reinhold, New
York, 1992. \item[{[DG]}] M. Demazure, and P. Gabriel, Groupes
alg\'{e}briques. Tome I: G\'{e}om\'{e}trie alg\'{e}brique,
g\'{e}n\'{e}ralit\'{e}s, groupes commutatifs, Masson, Paris;
North-Holland, Amsterdam, 1970. \item[{[EH]}] D. Eisenbud and J.
Harris, The geometry of schemes. Graduate Texts in Mathematics,
197. Springer-Verlag, New York, 2000. \item[{[GW]}] K. R. Goodearl
and R. B. Warfield, Jr.  An Introduction to Noncommutative
Noetherian Rings, London Math. Society Student Texts 16, Cambridge
University Press, 1989. \item[{[H]}] J.E. Humphreys, Linear
algebraic groups. Graduate Texts in Mathematics, 21.
Springer-Verlag, New York-Heidelberg, 1975. \item[{[Ja1]}] J. C.
Jantzen, Moduln mit einem h\"{o}chsten Gewicht. Lecture Notes in
Mathematics, 750. Springer, Berlin, 1979. \item[{[Ja2]}]
$\underline{\quad\quad\quad\quad\quad}$, Einh\"{u}llende Algebren
halbeinfacher Lie-Algebren.  Ergebnisse der Mathematik und ihrer
Grenzgebiete (3) Springer-Verlag, Berlin, 1983. \item[{[Ja3]}]
$\underline{\quad\quad\quad\quad\quad}$, Representations of
algebraic groups. Pure and Applied Mathematics, 131. Academic
Press, Inc., Boston, MA, 1987.
 \item[{[J1]}] A. Joseph,
On the associated variety of a primitive ideal, J. of Algebra 93 (1985),
509-523.
\item[{[J2]}] $\underline{\quad\quad\quad\quad\quad}$,  A surjectivity theorem for rigid highest weight
modules. Invent. Math.   92 (1988), no. 3, 567--596.
\item[{[J3]}] $\underline{\quad\quad\quad\quad\quad}$, Orbital varieties, Goldie rank
polynomials and unitary highest weight modules.  Algebraic and
analytic methods in representation theory (S\o nderborg, 1994),
53-98, Perspect. Math., 17, Academic Press, San Diego, CA, 1997.
\item[{[K]}] V. G. Kac, Lie superalgebras, Adv. in Math. 26 (1977), 8-96.
\item[{[KK]}] E. Kirkman, and J. Kuzmanovich, Minimal prime ideals in enveloping algebras of
Lie superalgebras.  Proc. Amer. Math. Soc.  124  (1996),
1693--1702.
\item[{[KL]}] G. R. Krause and T. H. Lenagan, Growth of algebras and
Gelfand-Kirillov dimension, Graduate Studies in Mathematics, vol.
22, American Mathematical Society, 2000.
\item[{[L]}] T.Y. Lam, The Algebraic Theory of quadratic Forms, Benjamin-Cummings,
Reading, Mass. 1980.
\item[{[Le]}]  D.A. Leites, Introduction to the theory
of supermanifolds. (Russian)  Uspekhi Mat. Nauk  35  (1980), 3-57,
Russian Math. Surveys 35 (1980), 1-64.
\item[{[LS]}] T. Levasseur, and S.P. Smith, Primitive ideals and nilpotent orbits in type $G\sb 2$,
 J. Algebra  114 (1988),  81--105, Corrigendum: J. Algebra  118  (1988), 261.
\item[{[Ma]}] Y. I. Manin, Gauge field theory and complex
geometry.  Second edition.  Grundlehren der Mathematischen
Wissenschaften, 289.  Springer-Verlag, Berlin, 1997.
\item[{[McR]}] J.C. McConnell, and J. C. Robson, Noncommutative
Noetherian rings. Graduate Studies in Mathematics, 30. American
Mathematical Society, Providence, RI, 2001. \item[{[Me]}] A.
Melnikov, Irreducibility of the associated varieties of simple
highest weight modules in $\FRAK{s}\FRAK{l}(n)$, C.R. Acad. Sci.
Paris Sr. I Math. 316 (1993), 53-57. \item[{[Mo]}] S. Montgomery,
Hopf algebras and their actions on rings. CBMS Regional Conference
Series in Mathematics, 82, American Mathematical Society,
Providence, RI, 1993. \item[{[M1]}] I. M. Musson, On the center of
the enveloping algebra of a classical simple Lie superalgebra, J.
Algebra 193 (1997), 75-101. \item[{[M2]}]
$\underline{\quad\quad\quad\quad\quad\quad}$, Associated varieties
for classical Lie superalgebras, pages 177-188 in Hopf algebras
and quantum groups (Brussels, 1998), 177--188, Lecture Notes in
Pure and Appl. Math., 209, Dekker, New York, 2000.
\item[{[M3]}]$\underline{\quad\quad\quad\quad\quad\quad}$, I.M.
Musson, Some Lie superalgebras associated to the Weyl algebras.
Proc. Amer. Math. Soc.  127  (1999),   2821--2827. \item[{[P]}] G.
Pinczon, The enveloping algebra of the Lie superalgebra ${\rm
osp}(1,2)$.  J. Algebra  132  (1990), 219--242. \item[{[Sch]}] M.
Scheunert, The Theory of Lie Superalgebras, Lecture Notes in
Mathematics, 716, Springer-Verlag, Berlin, 1979. \item[{[S1]}] V.
Serganova, On representations of the Lie superalgebra $p(n)$. J.
Algebra  258 (2002),   615--630. \item[{[S2]}]
$\underline{\quad\quad\quad\quad\quad\quad}$, A reduction method
for atypical representations of classical Lie superalgebras,
Advances in Mathematics 180, (2003), 248-274. \item[{[T]}] T.
Tanisaki, Characteristic varieties of highest weight modules and
primitive quotients. Representations of Lie groups, Kyoto,
Hiroshima, 1986, 1--30, Adv. Stud. Pure Math., 14, Academic Press,
Boston, MA, 1988. \item[{[Z]}]Y.M. Zou, Finite-dimensional
representations of $\Gamma(\sigma\sb 1,\sigma\sb 2,\sigma\sb 3)$.
J. Algebra  169  (1994),  no. 3, 827--846.
\end{itemize}
\end{document}